\documentclass[12pt]{article}

\catcode`\@=11

\thicklines
\newskip\Einheit \Einheit=.6cm
\newcount\xcoord \newcount\ycoord
\newdimen\xdim \newdimen\ydim \newdimen\PfadD@cke \newdimen\Pfadd@cke
\def\PfadDicke#1{\PfadD@cke#1 \divide\PfadD@cke by2 
\Pfadd@cke\PfadD@cke \multiply\PfadD@cke by2}
\long\def\LOOP#1\REPEAT{\def\BODY{#1}\ITERATE}
\def\ITERATE{\BODY \let\next\ITERATE \else\let\next\relax\fi \next}
\let\REPEAT=\fi
\def\Punkt{\hbox{\raise-2pt\hbox to0pt{\hss\scriptsize$\bullet$\hss}}}

\def\DuennPunkt(#1,#2){\unskip
  \raise#2 \Einheit\hbox to0pt{\hskip#1 \Einheit
          \raise-1.5pt\hbox to0pt{\hss\tiny$\bullet$\hss}\hss}}
		  
\def\NormalPunkt(#1,#2){\unskip
  \raise#2 \Einheit\hbox to0pt{\hskip#1 \Einheit
          \raise-3pt\hbox to0pt{\hss\large$\bullet$\hss}\hss}}
\def\DickPunkt(#1,#2){\unskip
  \raise#2 \Einheit\hbox to0pt{\hskip#1 \Einheit
          \raise-4pt\hbox to0pt{\hss\Large$\bullet$\hss}\hss}}
\def\Kreis(#1,#2){\unskip
  \raise#2 \Einheit\hbox to0pt{\hskip#1 \Einheit
          \raise-4pt\hbox to0pt{\hss\Large$\circ$\hss}\hss}}
\def\Diagonale(#1,#2)#3{\unskip\leavevmode
  \xcoord#1\relax \ycoord#2\relax
      \raise\ycoord \Einheit\hbox to0pt{\hskip\xcoord \Einheit
         \unitlength\Einheit
         \line(1,1){#3}\hss}}
\def\AntiDiagonale(#1,#2)#3{\unskip\leavevmode
  \xcoord#1\relax \ycoord#2\relax \advance\xcoord by -0.05\relax
      \raise\ycoord \Einheit\hbox to0pt{\hskip\xcoord \Einheit
         \unitlength\Einheit
         \line(1,-1){#3}\hss}}
\def\Pfad(#1,#2),#3\endPfad{\unskip\leavevmode
  \xcoord#1 \ycoord#2 \thicklines\ZeichnePfad#3\endPfad\thinlines}
\def\ZeichnePfad#1{\ifx#1\endPfad\let\next\relax
  \else\let\next\ZeichnePfad
    \ifnum#1=1
      \raise\ycoord \Einheit\hbox to0pt{\hskip\xcoord \Einheit
         \vrule height\Pfadd@cke width1 \Einheit depth\Pfadd@cke\hss}%
      \advance\xcoord by 1
     \else\ifnum#1=2
      \raise\ycoord \Einheit\hbox to0pt{\hskip\xcoord \Einheit
         \unitlength\Einheit
         \line(0,1){1}\hss}
      \advance\xcoord by 0
      \advance\ycoord by 1
 \else\ifnum#1=3
      \raise\ycoord \Einheit\hbox to0pt{\hskip\xcoord \Einheit
         \unitlength\Einheit
         \line(1,1){1}\hss}
      \advance\xcoord by 1
      \advance\ycoord by 1
    \else\ifnum#1=4
      \raise\ycoord \Einheit\hbox to0pt{\hskip\xcoord \Einheit
         \unitlength\Einheit
         \line(1,-1){1}\hss}
      \advance\xcoord by 1
      \advance\ycoord by -1
   \else\ifnum#1=5
      \raise\ycoord \Einheit\hbox to0pt{\hskip\xcoord \Einheit
         \unitlength\Einheit
         \line(2,1){2}\hss}
      \advance\xcoord by 2
      \advance\ycoord by 1
	  \else\ifnum#1=6
      \raise\ycoord \Einheit\hbox to0pt{\hskip\xcoord \Einheit
         \unitlength\Einheit
         \line(2,-1){2}\hss}
      \advance\xcoord by 2
      \advance\ycoord by -1
	  \else\ifnum#1=7
      \raise\ycoord \Einheit\hbox to0pt{\hskip\xcoord \Einheit
         \unitlength\Einheit
         \line(3,1){3}\hss}
      \advance\xcoord by 3
      \advance\ycoord by 1
	  \else\ifnum#1=8
      \raise\ycoord \Einheit\hbox to0pt{\hskip\xcoord \Einheit
         \unitlength\Einheit
         \line(3,-1){3}\hss}
      \advance\xcoord by 3
      \advance\ycoord by -1
    \fi\fi\fi\fi\fi\fi\fi\fi
  \fi\next}
\def\hSSchritt{\leavevmode\raise-.4pt\hbox 
to0pt{\hss.\hss}\hskip.2\Einheit
  \raise-.4pt\hbox to0pt{\hss.\hss}\hskip.2\Einheit
  \raise-.4pt\hbox to0pt{\hss.\hss}\hskip.2\Einheit
  \raise-.4pt\hbox to0pt{\hss.\hss}\hskip.2\Einheit
  \raise-.4pt\hbox to0pt{\hss.\hss}\hskip.2\Einheit}
\def\vSSchritt{\vbox{\baselineskip.2\Einheit\lineskiplimit0pt
\hbox{.}\hbox{.}\hbox{.}\hbox{.}\hbox{.}}}
\def\DSSchritt{\leavevmode\raise-.4pt\hbox to0pt{%
  \hbox to0pt{\hss.\hss}\hskip.2\Einheit
  \raise.2\Einheit\hbox to0pt{\hss.\hss}\hskip.2\Einheit
  \raise.4\Einheit\hbox to0pt{\hss.\hss}\hskip.2\Einheit
  \raise.6\Einheit\hbox to0pt{\hss.\hss}\hskip.2\Einheit
  \raise.8\Einheit\hbox to0pt{\hss.\hss}\hss}}
\def\dSSchritt{\leavevmode\raise-.4pt\hbox to0pt{%
  \hbox to0pt{\hss.\hss}\hskip.2\Einheit
  \raise-.2\Einheit\hbox to0pt{\hss.\hss}\hskip.2\Einheit
  \raise-.4\Einheit\hbox to0pt{\hss.\hss}\hskip.2\Einheit
  \raise-.6\Einheit\hbox to0pt{\hss.\hss}\hskip.2\Einheit
  \raise-.8\Einheit\hbox to0pt{\hss.\hss}\hss}}
\def\SPfad(#1,#2),#3\endSPfad{\unskip\leavevmode
  \xcoord#1 \ycoord#2 \ZeichneSPfad#3\endSPfad}
\def\ZeichneSPfad#1{\ifx#1\endSPfad\let\next\relax
  \else\let\next\ZeichneSPfad
    \ifnum#1=1
      \raise\ycoord \Einheit\hbox to0pt{\hskip\xcoord \Einheit
         \hSSchritt\hss}%
      \advance\xcoord by 1
    \else\ifnum#1=2
      \raise\ycoord \Einheit\hbox to0pt{\hskip\xcoord \Einheit
        \hbox{\hskip-2pt \vSSchritt}\hss}%
      \advance\ycoord by 1
    \else\ifnum#1=3
      \raise\ycoord \Einheit\hbox to0pt{\hskip\xcoord \Einheit
         \DSSchritt\hss}
      \advance\xcoord by 1
      \advance\ycoord by 1
    \else\ifnum#1=4
      \raise\ycoord \Einheit\hbox to0pt{\hskip\xcoord \Einheit
         \dSSchritt\hss}
      \advance\xcoord by 1
      \advance\ycoord by -1
    \fi\fi\fi\fi
  \fi\next}
\def\Koordinatenachsen(#1,#2){\unskip
 \hbox to0pt{\hskip-.5pt\vrule height#2 \Einheit width.5pt depth1 
\Einheit}%
 \hbox to0pt{\hskip-1 \Einheit \xcoord#1 \advance\xcoord by1
    \vrule height0.25pt width\xcoord \Einheit depth0.25pt\hss}}
\def\Koordinatenachsen(#1,#2)(#3,#4){\unskip
 \hbox to0pt{\hskip-.5pt \ycoord-#4 \advance\ycoord by1
    \vrule height#2 \Einheit width.5pt depth\ycoord \Einheit}%
 \hbox to0pt{\hskip-1 \Einheit \hskip#3\Einheit 
    \xcoord#1 \advance\xcoord by1 \advance\xcoord by-#3 
    \vrule height0.25pt width\xcoord \Einheit depth0.25pt\hss}}
\def\Gitter(#1,#2){\unskip \xcoord0 \ycoord0 \leavevmode
  \LOOP\ifnum\ycoord<#2
    \loop\ifnum\xcoord<#1
      \raise\ycoord \Einheit\hbox to0pt{\hskip\xcoord 
\Einheit\Punkt\hss}%
      \advance\xcoord by1
    \repeat
    \xcoord0
    \advance\ycoord by1
  \REPEAT}
\def\Gitter(#1,#2)(#3,#4){\unskip \xcoord#3 \ycoord#4 \leavevmode
  \LOOP\ifnum\ycoord<#2
    \loop\ifnum\xcoord<#1
      \raise\ycoord \Einheit\hbox to0pt{\hskip\xcoord 
\Einheit\Punkt\hss}%
      \advance\xcoord by1
    \repeat
    \xcoord#3
    \advance\ycoord by1
  \REPEAT}
\def\Label#1#2(#3,#4){\unskip \xdim#3 \Einheit \ydim#4 \Einheit
  \def\lo{\advance\xdim by-.5 \Einheit \advance\ydim by.5 \Einheit}%
  \def\llo{\advance\xdim by-.25cm \advance\ydim by.5 \Einheit}%
  \def\loo{\advance\xdim by-.5 \Einheit \advance\ydim by.25cm}%
  \def\o{\advance\ydim by.25cm}%
  \def\ro{\advance\xdim by.5 \Einheit \advance\ydim by.5 \Einheit}%
  \def\rro{\advance\xdim by.25cm \advance\ydim by.5 \Einheit}%
  \def\roo{\advance\xdim by.5 \Einheit \advance\ydim by.25cm}%
  \def\l{\advance\xdim by-.30cm}%
  \def\r{\advance\xdim by.30cm}%
  \def\lu{\advance\xdim by-.5 \Einheit \advance\ydim by-.6 \Einheit}%
  \def\llu{\advance\xdim by-.25cm \advance\ydim by-.6 \Einheit}%
  \def\luu{\advance\xdim by-.5 \Einheit \advance\ydim by-.30cm}%
  \def\u{\advance\ydim by-.30cm}%
  \def\ru{\advance\xdim by.5 \Einheit \advance\ydim by-.6 \Einheit}%
  \def\rru{\advance\xdim by.25cm \advance\ydim by-.6 \Einheit}%
  \def\ruu{\advance\xdim by.5 \Einheit \advance\ydim by-.30cm}%
  #1\raise\ydim\hbox to0pt{\hskip\xdim
     \vbox to0pt{\vss\hbox to0pt{\hss$#2$\hss}\vss}\hss}%
}
\catcode`\@=12

\usepackage{amsmath,amsthm}
\usepackage{amssymb}
\usepackage{color}
\usepackage{xspace}
\usepackage{pstricks,pst-node}
\usepackage[colorlinks=true,
linkcolor=green,
filecolor=brown,
citecolor=green]{hyperref}
\definecolor{gray}{rgb}{.7,.7,.7}

\def\blue1{\textcolor{blue} }

\def\gray{\textcolor{gray} }
\def\IOTs{increasing ordered trees\xspace}
\def\IOT{increasing ordered tree\xspace}
\def\HL{height-labeled\xspace}
\def\v{\vert}
\def\si{\sigma}
\def\gf{generating function\xspace}
\def\a{\ensuremath{\mathcal A}\xspace}
\def\b{\ensuremath{\mathcal B}\xspace}
\def\C{\ensuremath{\mathcal C}\xspace}

\def\p{\ensuremath{\mathcal P}\xspace}
\def\u{\ensuremath{\mathcal U}\xspace}
\def\V{\ensuremath{\mathcal V}\xspace}
\def\si{\sigma}

\def\mbf#1{\mathchoice{\hbox{\boldmath $\displaystyle #1$}}
        {\hbox{\boldmath $\textstyle #1$}}
        {\hbox{\boldmath $\scriptstyle #1$}}
        {\hbox{\boldmath $\scriptscriptstyle #1$}}} 

\newcommand{\StirlingCycle}[2]{\genfrac{[}{]}{0pt}{}{#1}{#2}}
\newcommand{\SecondEulerian}[2]{\left\langle\!\!\! 
     \genfrac{\langle}{\rangle}{0pt}{}{#1}{#2} \!\!\!\right\rangle}
\newcommand{\SecondEulerianInline}[2]{\left\langle \!\!
     \genfrac{\langle}{\rangle}{0pt}{}{#1}{#2}\!\! \right\rangle}

\begin{document}
\newtheorem{theorem}{Theorem}
\newtheorem{defn}[theorem]{Definition}
\newtheorem{lemma}[theorem]{Lemma}
\newtheorem{prop}[theorem]{Proposition}
\newtheorem{cor}[theorem]{Corollary}

\begin{center}
{\Large
A Combinatorial Survey of Identities for the Double Factorial      \\ 
}

\vspace{10mm}
DAVID CALLAN  \\
Department of Statistics  \\
\vspace*{-1mm}
University of Wisconsin-Madison  \\
\vspace*{-1mm}
Medical Science Center \\
\vspace*{-1mm}
1300 University Ave  \\
\vspace*{-1mm}
Madison, WI \ 53706-1532  \\
{\bf callan@stat.wisc.edu}  \\
\vspace{5mm}

June 6, 2009
\end{center}

\begin{abstract}
We survey combinatorial interpretations of some dozen identities for the double factorial 
such as, for instance,
$(2n-2)!! + \sum_{k=2}^{n} \frac{(2n-1)!!(2k-4)!!}{(2k-1)!!} = 
(2n-1)!!$. Our methods are mostly bijective.
\end{abstract}

\section{Introduction}

There are a surprisingly large number of identities for the odd double factorial 
$(2n-1)!!=(2n-1)\cdot (2n-3)\cdots 3 \cdot 1 =\frac{(2n)!}{2^{n}n!} $ that involve round numbers (small 
prime factors), as well as several that don't.
The purpose of this paper is to present (and in some cases, review) combinatorial interpretations 
of these identities. Section \ref{manifestation} reviews combinatorial 
constructs counted by $(2n-1)!!$. Section \ref{Hafnian} uses Hafnians 
to establish one of these manifestations. The subsequent sections contain 
the main  results and treat 
individual identities, presenting one or more combinatorial interpretations 
for each; Section \ref{round} is devoted to 
round-number identities, Section \ref{other} to non-round identities, 
Section \ref{refine} to refinements involving double summations 
interpreted by two statistics in addition to size.

The even double factorial is
$(2n)!!=2n\cdot (2n-2) \cdots 2=2^{n}n!.$ The recurrence $k!!=k(k-2)!!$ allows the definition of the 
double factorial to be extended to odd negative arguments. In 
particular, the values $(-1)!!=1$ and $(-3)!!=-1$ will arise in 
some of the identities. We use the notations $[n]$ for the set 
$\{1,2,\ldots,n\}$ and $[a,b]$ for the closed interval of integers 
from $a$ to $b$. For negative $k$, $\binom{n}{k}=0$ as usual, except that 
identities (\ref{rightPathFromRoot}) and (\ref{firstascent}) below require 
$\binom{-1}{-1}:=1$. It has become customary  
to draw trees down 
but when a construction involves growing a tree, it seems more natural 
to draw it up. Also, to visualize ordered trees as Dyck paths, they 
must go up. So we will combine arborological pictures conjuring roots 
and leaves with the usual genealogical terminology of parents, 
children, and siblings. We sometimes refer to (clockwise) 
\emph{walkaround} order of the edges/vertices in a tree; more formally, it is the order 
edges are visited in depth-first search and the preorder of the 
vertices.

\section[Combinatorial manifestations of (2$n-1$)!!]{Combinatorial 
manifestations of (2$\protect\mbf{n}-\protect\mbf{1}$)!!}
\label{manifestation}

We begin with a review of some combinatorial manifestations of $(2n-1)!!$, 
each illustrated for the case $n=2$. In all cases, the parameter $n$ 
is the \emph{size} of the object.
\subsection{Trapezoidal words} \label{symTrap} 
\vspace*{-5mm}
The Cartesian product
$[1]\times [3]\times \ldots \times [2n-1]$. 
\[
11,\quad 12,\quad 13.
    \]
The elements of this Cartesian product, the most obvious manifestation of $(2n-1)!!$, 
were called \emph{trapezoidal} words by Riordan \cite{rior76}. A minor 
variation is \emph{symmetric trapezoidal} words: the Cartesian product
$[0,0]\times [-1,1]\times [-2,2]\times \ldots \times [-(n-1),n-1)]$.

\subsection{Perfect matchings}  
\vspace*{-5mm}
Perfect matchings of $[2n]$.
    \[
    12/34,\quad 13/24,\quad 14/23.
    \]
A \emph{perfect matching} of $[2n]=\{1,2,\ldots,2n\}$ is a partition of $[2n]$ 
into 2-element subsets $a(1)\,b(1)\, / \, a(2)\,b(2)\, / \, \ldots   
\, / \, a(n)\, b(n)$ written, in standard form, so that $a(i)<b(i)$ 
for all $i$, and 
$a(1)<a(2)<\ldots <a(n)$. Erasing the virgules (slashes) then gives a 
bijection to the \emph{perfect matching permutations} of $[2n]$, denoted 
$\p_{n}$: the permutations $a(1)\,b(1)\,a(2)\,b(2)\, \ldots\, a(n)\,b(n)$ 
(denoted $(\mathbf{a},\mathbf{b})$) 
of $[2n]$ satisfying $a(i)<b(i)$ for all $i$, and 
$a(1)<a(2)<\ldots <a(n)$. Given a perfect matching permutation 
$(\mathbf{a},\mathbf{b})$, form a list whose $a(i)$-th and $b(i)$-th 
entries are both $i$. This is a bijection to the permutations of the 
multiset  $\{1,1,2,2,\ldots,n,n\}$ in                     
which the first occurrences of $1,2,\ldots,n$ occur in that order.

A  perfect matching of $[2n]$ can be 
regarded as a fixed-point-free involution on $[2n]$ and also as a 1-regular 
graph on $[2n]$, whose pictorial representation is sometimes called a 
Brauer diagram.
\subsection{Stirling permutations}
\vspace*{-5mm}
Permutations of the multiset $\{1,1,2,2,\ldots,n,n\}$ in 
    which, for each $i$, all entries between the two occurrences of 
    $i$ exceed $i$ \cite{stirpoly}.
    \[
    1122,\quad 1221,\quad 2211.
     \] 

\subsection{Increasing ordered trees} \label{incrTree}
\vspace*{-5mm}
Increasing vertex-labeled ordered trees of $n$ edges, 
    label set $[0,n]$ \cite{klazar1,klazar2}.  
\Einheit=0.5cm
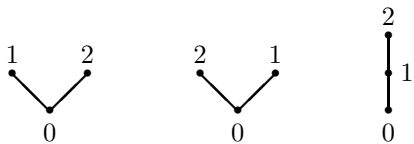
\[
\Pfad(-6,1),43\endPfad
\Pfad(-1,1),43\endPfad
\Pfad(4,0),22\endPfad
\DuennPunkt(-6,1)
\DuennPunkt(-5,0)
\DuennPunkt(-4,1)
\DuennPunkt(-1,1)
\DuennPunkt(0,0)
\DuennPunkt(1,1)
\DuennPunkt(4,0)
\DuennPunkt(4,1)
\DuennPunkt(4,2)
\Label\o{ \textrm{\footnotesize 1}}(-6,1)
\Label\o{ \textrm{\footnotesize 2}}(-4,1)
\Label\o{ \textrm{\footnotesize 2}}(-1,1)
\Label\o{ \textrm{\footnotesize 1}}(1,1)
\Label\r{ \textrm{\footnotesize 1}}(3.9,1)
\Label\o{ \textrm{\footnotesize 2}}(4,2)
\Label\u{ \textrm{\footnotesize 0}}(-5,0)
\Label\u{ \textrm{\footnotesize 0}}(0,0)
\Label\u{ \textrm{\footnotesize 0}}(4,0)
\]
Given such a tree, delete the root label and transfer the remaining labels from vertices to 
parent edges. Walk clockwise around the tree thereby traversing each 
edge twice and record labels in the order encountered. This is a 
bijection to Stirling permutations due to Svante Janson \cite{janson08}. 
Label sets other than $[0,n]$ may arise and the term \emph{standard} then 
emphasizes that a tree's label set is an initial segment of the 
nonnegative integers.  
\subsection{Leaf-labeled 0-2 trees} 
\vspace*{-5mm}
0-2 trees (rooted, unordered, each vertex has 0 or 2 
    children) with $n+1$ labeled leaves, label set $[1,n+1]$ 
    \cite[Chapter 5.2.6]{ec2}.
    
\Einheit=0.5cm
\[
\Pfad(-7,2),43\endPfad
\Pfad(-6,1),43\endPfad
\Pfad(-2,2),43\endPfad
\Pfad(-1,1),43\endPfad
\Pfad(3,2),43\endPfad
\Pfad(4,1),43\endPfad
\DuennPunkt(-7,2)
\DuennPunkt(-6,1)
\DuennPunkt(-5,0)
\DuennPunkt(-5,2)
\DuennPunkt(-4,1)
\DuennPunkt(-2,2)
\DuennPunkt(-1,1)
\DuennPunkt(0,0)
\DuennPunkt(0,2)
\DuennPunkt(1,1)
\DuennPunkt(3,2)
\DuennPunkt(4,1)
\DuennPunkt(5,0)
\DuennPunkt(5,2)
\DuennPunkt(6,1)
\Label\o{ \textrm{\footnotesize 1}}(-7,2)
\Label\o{ \textrm{\footnotesize 2}}(-5,2)
\Label\o{ \textrm{\footnotesize 3}}(-4,1)
\Label\o{ \textrm{\footnotesize 1}}(-2,2)
\Label\o{ \textrm{\footnotesize 3}}(0,2)
\Label\o{ \textrm{\footnotesize 2}}(1,1)
\Label\o{ \textrm{\footnotesize 2}}(3,2)
\Label\o{ \textrm{\footnotesize 3}}(5,2)
\Label\o{ \textrm{\footnotesize 1}}(6,1)
\]

\subsection{Height-labeled Dyck paths}  \label{htdyck} 
\vspace*{-5mm}
Dyck paths \cite[Exercise 6.19]{ec2} of 
    $n$ upsteps and $n$ downsteps in which each upstep is labeled 
    with a positive integer $\le $ the height of its top vertex.
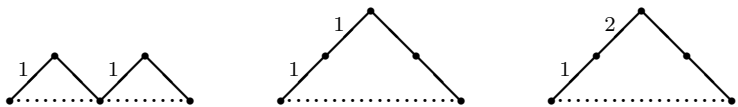
\Einheit=0.6cm
\[
\Pfad(-8,0),3434\endPfad
\SPfad(-8,0),1111\endSPfad
\Pfad(-2,0),3344\endPfad
\SPfad(-2,0),1111\endSPfad
\Pfad(4,0),3344\endPfad
\SPfad(4,0),1111\endSPfad
\DuennPunkt(-8,0)
\DuennPunkt(-7,1)
\DuennPunkt(-6,0)
\DuennPunkt(-5,1)
\DuennPunkt(-4,0)
\DuennPunkt(-2,0)
\DuennPunkt(-1,1)
\DuennPunkt(0,2)
\DuennPunkt(1,1)
\DuennPunkt(2,0)
\DuennPunkt(4,0)
\DuennPunkt(5,1)
\DuennPunkt(6,2)
\DuennPunkt(7,1)
\DuennPunkt(8,0)
\Label\o{\textrm{{\scriptsize 1}}}(-7.7,0.3)
\Label\o{\textrm{{\scriptsize 1}}}(-5.7,0.3)
\Label\o{\textrm{{\scriptsize 1}}}(-1.7,0.3)
\Label\o{\textrm{{\scriptsize 1}}}(-0.7,1.3)
\Label\o{\textrm{{\scriptsize 1}}}(4.3,0.3)
\Label\o{\textrm{{\scriptsize 2}}}(5.3,1.3)
\]
This item is due to Jean Fran\c{c}on and G\'{e}rard Viennot 
\cite{francon78,francon79}; they observe that it is a consequence of 
their bijection from permutations to certain marked-up lattice paths. 
Here is perhaps the simplest proof and several further bijections that 
prove the result will appear in the course of the paper.  
A height-labeled (HL) Dyck path $P$ of size $n-1$ has $2n-1$ 
vertices, each of which can be used to construct a height-labeled Dyck path of 
size $n$: split $P$ at the specified vertex into subpaths $P_{1}$ and 
$P_{2}$, insert an upstep between $P_{1}$ and $P_{2}$, increment by 1 
the labels on $P_{2}$, and append a downstep, as illustrated below.
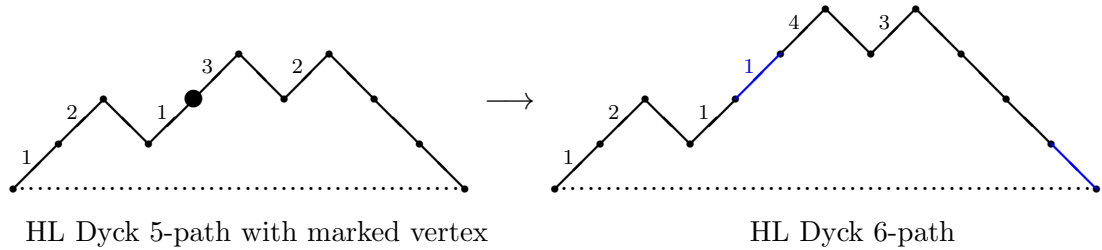
\Einheit=0.6cm
\[
\Pfad(-12,0),3343343444\endPfad
\Pfad(0,0),3343\endPfad
\Pfad(5,3),343444\endPfad
\SPfad(-12,0),1111111111\endSPfad
\SPfad(0,0),111111111111\endSPfad
\DuennPunkt(-12,0)
\DuennPunkt(-11,1)
\DuennPunkt(-10,2)
\DuennPunkt(-9,1)
\DickPunkt(-8,2)
\DuennPunkt(-7,3)
\DuennPunkt(-6,2)
\DuennPunkt(-5,3)
\DuennPunkt(-4,2)
\DuennPunkt(-3,1)
\DuennPunkt(-2,0)
\DuennPunkt(0,0)
\DuennPunkt(1,1)
\DuennPunkt(2,2)
\DuennPunkt(3,1)
\DuennPunkt(4,2)
\DuennPunkt(5,3)
\DuennPunkt(6,4)
\DuennPunkt(7,3)
\DuennPunkt(8,4)
\DuennPunkt(9,3)
\DuennPunkt(10,2)
\DuennPunkt(11,1)
\DuennPunkt(12,0)
\Label\o{\longrightarrow}(-1,1.5)
\Label\o{\textrm{{\scriptsize 1}}}(-11.7,0.3)
\Label\o{\textrm{{\scriptsize 2}}}(-10.7,1.3)
\Label\o{\textrm{{\scriptsize 1}}}(-8.7,1.3)
\Label\o{\textrm{{\scriptsize 3}}}(-7.7,2.3)
\Label\o{\textrm{{\scriptsize 2}}}(-5.7,2.3)
\Label\o{\textrm{{\scriptsize 1}}}(0.3,0.3)
\Label\o{\textrm{{\scriptsize 2}}}(1.3,1.3)
\Label\o{\textrm{{\scriptsize 1}}}(3.3,1.3)
\Label\o{\textrm{{\scriptsize 4}}}(5.3,3.3)
\Label\o{\textrm{{\scriptsize 3}}}(7.3,3.3)
\blue1{\Pfad(4,2),3\endPfad
\Pfad(11,1),4\endPfad
\Label\o{\textrm{{\scriptsize 1}}}(4.3,2.3) }
\Label\o{\textrm{{\small HL Dyck 5-path with marked 
vertex}}}(-6.6,-1.4)
\Label\o{\textrm{{\small HL Dyck 6-path}}}(6.6,-1.4)
\]

\vspace*{1mm}

\noindent The process can be reversed by locating the last upstep with label 1 
in a size-$n$ path. Thus there is a multiplying factor of $2n-1$ from 
size $n-1$ to size $n$ and the number of height-labeled Dyck $n$-paths is indeed 
$(2n-1)!!$.

\subsection{Height-labeled ordered trees}  \label{htTree} 
\vspace*{-5mm}
Ordered trees of 
    $n$ edges in which each non-root vertex is labeled 
    with a positive integer $\le $ its height (distance from root).
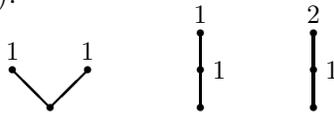
\Einheit=0.5cm
\[
\Pfad(-4,1),43\endPfad
\Pfad(1,0),22\endPfad
\Pfad(4,0),22\endPfad
\DuennPunkt(-4,1)
\DuennPunkt(-3,0)
\DuennPunkt(-2,1)
\DuennPunkt(1,0)
\DuennPunkt(1,1)
\DuennPunkt(1,2)
\DuennPunkt(4,0)
\DuennPunkt(4,1)
\DuennPunkt(4,2)
\Label\o{ \textrm{\footnotesize 1}}(-4,1)
\Label\o{ \textrm{\footnotesize 1}}(-2,1)
\Label\r{ \textrm{\footnotesize 1}}(0.9,1)
\Label\o{ \textrm{\footnotesize 1}}(1,2)
\Label\r{ \textrm{\footnotesize 1}}(3.9,1)
\Label\o{ \textrm{\footnotesize 2}}(4,2)
\]    
\noindent The ``accordion'' bijection from ordered trees to Dyck paths---burrow 
up the branches from the root and open out the tree as 
illustrated---sends non-root vertices to tops of upsteps and preserves 
height.
\begin{center} 
\begin{pspicture}(-8,-1)(16,1.5)
\psset{unit=.6cm}
\psline(-12,2)(-12,1)(-11,0)(-10,1)(-11,2)
\psline(-10,2)(-10,1)(-9,2)(-9,3)

\psdots(-12,2)(-12,1)(-11,0)(-10,1)(-11,2)(-10,2)(-9,2)(-9,3)

\psline(-6.1,0)(-7.1,1)(-7,2)(-6.9,1)(-6,.2)(-5.1,1)(-6,2)
(-5.1,1.2)(-5,2)(-4.9,1.2)(-4.1,2)(-4,3)(-3.9,2)(-4.9,1)(-5.9,0)

\rput(-8,1){$\longrightarrow$}
\rput(-3,1){$\longrightarrow$}

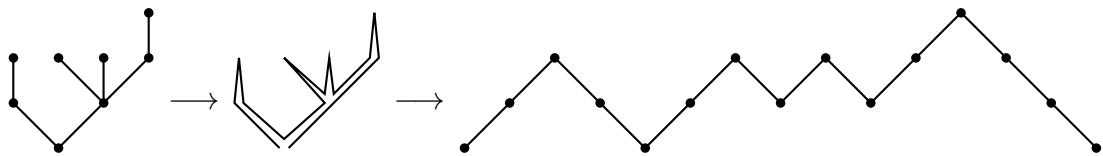
\rput(0,-1.5){{\small The ``accordion'' bijection}}
 
\psline(-2,0)(-1,1)(0,2)(1,1)(2,0)(3,1)(4,2)(5,1)(6,2)(7,1)(8,2)(9,3)(10,2)(11,1)(12,0)
\psdots(-2,0)(-1,1)(0,2)(1,1)(2,0)(3,1)(4,2)(5,1)(6,2)(7,1)(8,2)(9,3)(10,2)(11,1)(12,0)

\end{pspicture}
\end{center}  
Thus height-labeled ordered trees transparently correspond to height-labeled 
Dyck paths.

\subsection{Overhang paths}
\vspace*{-5mm}
Lattice paths of steps $(1,1),(1,-1),(-1,1)$ from (0,0) to $(2n,0)$ 
that lie in the first quadrant and do not self-intersect 
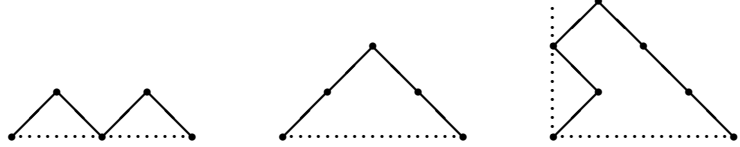
\cite{overhang}.
\Einheit=0.6cm
\[
\Pfad(-8,0),3434\endPfad
\SPfad(-8,0),1111\endSPfad
\Pfad(-2,0),3344\endPfad
\SPfad(-2,0),1111\endSPfad
\Pfad(4,0),3\endPfad
\Pfad(4,2),3444\endPfad
\Pfad(4,2),4\endPfad
\SPfad(4,0),1111\endSPfad
\SPfad(4,0),222\endSPfad
\DuennPunkt(-8,0)
\DuennPunkt(-7,1)
\DuennPunkt(-6,0)
\DuennPunkt(-5,1)
\DuennPunkt(-4,0)
\DuennPunkt(-2,0)
\DuennPunkt(-1,1)
\DuennPunkt(0,2)
\DuennPunkt(1,1)
\DuennPunkt(2,0)
\DuennPunkt(4,0)
\DuennPunkt(4,2)
\DuennPunkt(5,3)
\DuennPunkt(5,1)
\DuennPunkt(6,2)
\DuennPunkt(7,1)
\DuennPunkt(8,0)
\]
Listing the ordinates of the upstep tops 
is a bijection to trapezoidal words.

\section{Pfaffians, Hafnians and Dyck paths}\label{Hafnian}

The Pfaffian is usually defined for a skew-symmetric matrix but it 
can just as well be defined for the upper triangular array 
$T=(x_{ij})_{1\le i<j \le 2n}$\,:
\[
\textrm{Pf\,}(T)=\sum_{(\mathbf{a},\mathbf{b})\in \p_{n}} \textrm{sgn}
((\mathbf{a},\mathbf{b}))x_{a(1)b(1)}x_{a(2)b(2)}\cdots x_{a(n)b(n)},
\]
 a sum over all $(2n-1)!!$ perfect matching permutations 
$(\mathbf{a},\mathbf{b})$ in $\p_{n}$ (where $\textrm{sgn}((\mathbf{a},\mathbf{b}))$ is the sign of the 
permutation). 
The Hafnian of $T$ is given by the same sum but without the signs: 
\[
\textrm{Hf\,}(T)=\sum_{(\mathbf{a},\mathbf{b})\in \p_{n}} x_{a(1)b(1)}x_{a(2)b(2)}\cdots x_{a(n)b(n)}
\]
Obviously, the Hafnian of the all 1s array is $(2n-1)!!$.

\begin{prop}\label{hf}
  For an array   $T=(x_{ij})_{1\le i<j \le 2n}$ with constant rows 
  $x_{ij}:=x_{i}$, 
\[
\mbox{\emph{Hf\,}}(T)=\sum_{\mathbf{a} \in \a_{n}}\prod_{i=1}^{n}(2i-a(i))x_{a(i)}
\]
where $\a_{n}$ denotes the set of increasing sequences 
$\mathbf{a}=(a(1),a(2),\ldots,a(n))$ satisfying $1\le a(i) \le 2i-1$.
\end{prop}

Proof. The map ($\mathbf{a},\mathbf{b}) \rightarrow \mathbf{a}$ maps 
perfect matching permutations onto $\a_{n}$. Thus every term in 
$\textrm{Hf\,}(T)$ has the form $\prod_{i=1}^{n}x_{a(i)}$ for some 
$\mathbf{a} \in \a_{n}$ and the question is, how many of each form? Given 
$\mathbf{a} \in \a_{n}$, set $R=\{a(1),a(2),\ldots,a(n)\}$ and 
$C=[2n]\,\backslash\, R$, and let $C_{i}=\{c\in C\,:\,c>a(i)\}$ and 
$c_{i}=\v\,C_{i}\,\v$. Clearly, $C_{n}\subseteq C_{n-1}\subseteq 
\ldots \subseteq C_{1}$. The $b_{i}$'s in a perfect matching permutation
$(\mathbf{a},\mathbf{b})$ are subject only to the two restrictions: 
$b_{i}\in C_{i}$ and all $b_{i}$'s distinct. Thus there are $c_{n}$ 
choices for $b_{n},\ c_{n-1}-1$ choices for $b_{n-1},\ c_{n-2}-2$ 
choices for $b_{n-2}$ and so on. Hence the coefficient of 
$\prod_{i=1}^{n}x_{a(i)}$ in $\textrm{Hf\,}(T)$ is 
$\prod_{i=1}^{n}\big(c_{i}-(n-i)\big)$ and a simple check shows that 
$c_{i}-(n-i)=2i-a(i)$. \qed

A Dyck path can be coded by the positions in the path of its upsteps 
and this coding is a bijection from Dyck $n$-paths onto $\a_{n}$. A Dyck 
path can also be coded by the heights of the tops of its upsteps, 
giving  a bijection from Dyck $n$-paths to the sequences 
$\b_{n}=\{\big(b(i)\big)_{i=1}^{n}\}$ satisfying $b(1)=1$ and $1 \le b(i+1) \le 
b(i)+1$ for $1\le i \le n-1$. The two codes are related by the 
equality 
$b(i)=2i-a(i)$ for all $i$.

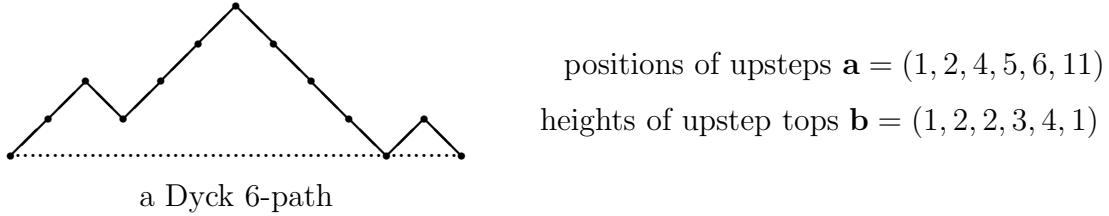
\vspace*{-4mm}

\Einheit=0.5cm
\[
\Pfad(-14,0),334333444434\endPfad
\SPfad(-14,0),111111111111\endSPfad
\DuennPunkt(-14,0)
\DuennPunkt(-13,1)
\DuennPunkt(-12,2)
\DuennPunkt(-11,1)
\DuennPunkt(-10,2)
\DuennPunkt(-9,3)
\DuennPunkt(-8,4)
\DuennPunkt(-7,3)
\DuennPunkt(-6,2)
\DuennPunkt(-5,1)
\DuennPunkt(-4,0)
\DuennPunkt(-3,1)
\DuennPunkt(-2,0)
\Label\u{ \textrm{positions of upsteps $\mathbf{a}=(1,2,4,5,6,11)$}}(7.9,3)
\Label\u{ \textrm{heights of upstep tops 
$\mathbf{b}=(1,2,2,3,4,1)$}}(7.5,1.5)
\Label\u{\textrm{a Dyck 6-path}}(-8,-0.5)
\]

\vspace*{2mm}

The correspondences, Dyck $n$-paths 
$\rightarrow \a_{n} \rightarrow \b_{n}$, along with the assertion of Prop. \ref{hf} and the fact that 
the Hafnian of the all 1s array is $(2n-1)!!$, now establish item 
\ref{htdyck}.

This proof can be worked up into a bijection from height-labeled Dyck 
paths to perfect matchings $(\mathbf{a},\mathbf{b})$. The positions 
of the upsteps in the Dyck path give $\mathbf{a}$, and $\mathbf{b}$ is 
formed from the positions of the downsteps, using the sequence of height labels 
$\mathbf{h}=\big(h(i)\big)_{i=1}^{n}$, as follows. First, write the downstep 
position list $\mathbf{d}$ in decreasing order. Then $b(n)$ is the $h(n)$-th 
entry in $\mathbf{d}$, $b(n-1)$ is the $h(n-1)$-th entry in the remaining $n-1$ 
elements of $\mathbf{d}$, $b(n-2)$ is the $h(n-2)$-th entry in the remaining $n-2$ 
elements and so on. An example is illustrated.

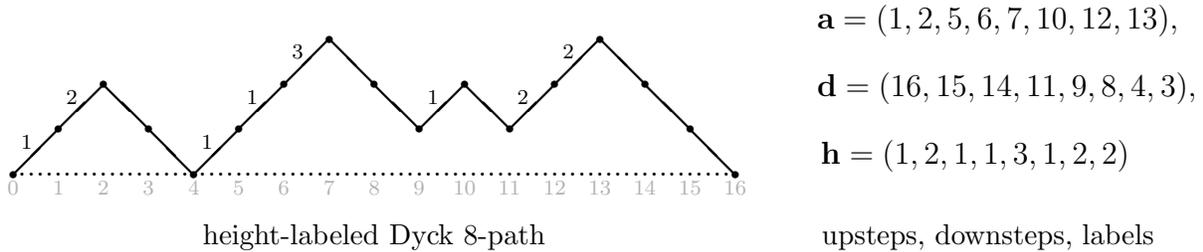
\Einheit=0.6cm
\[
\Pfad(-13,0),3344333443433444\endPfad
\SPfad(-13,0),1111111111111111\endSPfad
\DuennPunkt(-13,0)
\DuennPunkt(-12,1)
\DuennPunkt(-11,2)
\DuennPunkt(-10,1)
\DuennPunkt(-9,0)
\DuennPunkt(-8,1)
\DuennPunkt(-7,2)
\DuennPunkt(-6,3)
\DuennPunkt(-5,2)
\DuennPunkt(-4,1)
\DuennPunkt(-3,2)
\DuennPunkt(-2,1)
\DuennPunkt(-1,2)
\DuennPunkt(0,3)
\DuennPunkt(1,2)
\DuennPunkt(2,1)
\DuennPunkt(3,0)
\gray{\Label\u{\textrm{{\scriptsize 0}}}(-13,0.2)
\Label\u{\textrm{{\scriptsize 1}}}(-12,0.2)
\Label\u{\textrm{{\scriptsize 2}}}(-11,0.2)
\Label\u{\textrm{{\scriptsize 3}}}(-10,0.2)
\Label\u{\textrm{{\scriptsize 4}}}(-9,0.2)
\Label\u{\textrm{{\scriptsize 5}}}(-8,0.2)
\Label\u{\textrm{{\scriptsize 6}}}(-7,0.2)
\Label\u{\textrm{{\scriptsize 7}}}(-6,0.2)
\Label\u{\textrm{{\scriptsize 8}}}(-5,0.2)
\Label\u{\textrm{{\scriptsize 9}}}(-4,0.2)
\Label\u{\textrm{{\scriptsize 10}}}(-3,0.2)
\Label\u{\textrm{{\scriptsize 11}}}(-2,0.2)
\Label\u{\textrm{{\scriptsize 12}}}(-1,0.2)
\Label\u{\textrm{{\scriptsize 13}}}(0,0.2)
\Label\u{\textrm{{\scriptsize 14}}}(1,0.2)
\Label\u{\textrm{{\scriptsize 15}}}(2,0.2)
\Label\u{\textrm{{\scriptsize 16}}}(3,0.2) }
\Label\o{\textrm{{\scriptsize 1}}}(-12.7,0.3)
\Label\o{\textrm{{\scriptsize 2}}}(-11.7,1.3)
\Label\o{\textrm{{\scriptsize 1}}}(-8.7,0.3)
\Label\o{\textrm{{\scriptsize 1}}}(-7.7,1.3)
\Label\o{\textrm{{\scriptsize 3}}}(-6.7,2.3)
\Label\o{\textrm{{\scriptsize 1}}}(-3.7,1.3)
\Label\o{\textrm{{\scriptsize 2}}}(-1.7,1.3)
\Label\o{\textrm{{\scriptsize 2}}}(-0.7,2.3)
\Label\o{ \textrm{$\mathbf{a}=(1,2,5,6,7,10,12,13),$}}(8.8,3)
\Label\o{ \textrm{$\mathbf{d}=(16,15,14,11,9,8,4,3),$}}(9,1.5)
\Label\o{ \textrm{$\mathbf{h}=(1,2,1,1,3,1,2,2)$}}(8.3,0)
\Label\o{\textrm{{\small height-labeled Dyck 8-path}}}(-5,-1.8)
\Label\o{\textrm{{\small upsteps, downsteps, labels}}}(8.6,-1.8)
\]

\vspace*{3mm}

\noindent So $b(8)$ is the second entry of $\mathbf{d}$, namely 15; $b(7)$ is the second entry of
$(16,14,11,\ldots)$, namely 14;  $b(3)=16$ and so on. The result is 
$\mathbf{b}=(4,3,9,11,8,16,14,15)$ and the perfect matching is
\[
1\ 4\,/\,2\ 3\,/\,5\ 9\,/\,6\ 11\,/\,7\ 8\,/\,10\ 16\,/\,12\ 14\,/\,13\ 15.
\]

There is an analogue of Prop. \ref{hf} for the Pfaffian.
\begin{prop}\label{pf}
  For an array   $T=(x_{ij})_{1\le i<j \le 2n}$ with constant rows 
  $x_{ij}:=x_{i}$, 
\[
\mbox{\emph{Pf\,}}(T)=x_{1}x_{3}x_{5}\ldots x_{2n-1}.
\]
\end{prop}

Proof. For a perfect matching permutation ($\mathbf{a},\mathbf{b}$), take 
the smallest $i$ for which $a(i)$ and $b(i)$ are not consecutive 
integers (if there is one). Then the entries in 
($\mathbf{a},\mathbf{b}$) up through $a(i)$ 
necessarily form an initial segment of the positive integers ending at 
$a(i)=2i-1$, and $a(i+1)=2i$. So interchanging $b(i)$ and $b(i+1)$ 
gives another perfect matching permutation. Both contribute the same 
product to the Pfaffian but with opposite signs and hence they cancel out. 
The only surviving permutation under this involution is the identity, which contributes 
$x_{1}x_{3}x_{5}\ldots x_{2n-1}$. \qed

\begin{cor}
    For the array   $T=(x_{ij})_{1\le i<j \le 2n}$ with 
  $x_{ij}:=i$, 
\[
\mbox{\emph{Pf\,}}(T)=(2n-1)!!.
\]
\end{cor}
Proof. Put $x_{i}=i$ in Prop. \ref{pf}.

\section{Round-Number Identities}\label{round}
\subsection{}
\vspace*{-7mm}
\begin{equation}
 \sum_{k=0}^{n-1}\binom{n}{k+1} (2k-1)!!(2n-2k-3)!! = (2n-1)!!
    \label{leftSubtreeOfRoot}
\end{equation}
This identity counts increasing ordered trees of size $n$ by size  
$k$ of the leftmost subtree of the root. To see this, simply condition 
on the vertex set of the leftmost subtree. 
The bivariate \gf $\sum_{n\ge 1,k\ge 0}\binom{n}{k+1} (2k-1)!!(2n-2k-3)!! 
\frac{\textrm{{\,\normalsize $x$}}^{n}}{\textrm{{\small $n$}}!} y^{k}$ is 
\[
 \frac{1-\sqrt{1-2xy}}{y \sqrt{1-2x}},
\]
and the first few values are
\[
\begin{array}{c|ccccc}
	n^{\textstyle{\,\backslash \,k}} & 0 & 1 & 2 & 3 & 4  \\
\hline 
	1&    1 &   & & & \\
 	2&    2 & 1 & & &   \\
	3&    9 & 3 & 3 & &  \\
	4&    60 & 18  & 12 & 15 &   \\ 
	5&    525 & 150 & 90 & 75 & 105  \\
 
 \end{array}
\]

\subsection{}\label{rightPath}
\vspace*{-4mm}

\begin{equation}
 \sum_{k=0}^{n}\binom{2n-k-1}{k-1} 
 \frac{(2n-2k-1)(2n-k+1)}{k+1}(2n-2k-3)!!  = (2n-1)!!
    \label{rightPathFromRoot}
\end{equation}
This identity counts increasing ordered trees of size $n$ by length  
$k$ of the rightmost path from the root---the path that starts at the 
root and successively goes to the rightmost child until it reaches a 
leaf. To see this, let $u(n,k)$ 
be the number of such trees. Clearly, $u(0,0)=1$ and, for $n\ge 1,\ u(n,1)=n(2n-3)!!$ since there 
are $n$ choices for the rightmost child of the root and this child is 
a leaf. Now suppose $n\ge k\ge 2$ and consider 
trees of size $n-1$. If the rightmost path has 
length $\ge k-1$ then there is just one way to insert $n$ to get a 
size-$n$ tree with rightmost path of length $k$, and $n$ then ends 
the rightmost path. On the other hand, if the rightmost path has 
length $k$, then adding $n$ anywhere except as the rightmost child of 
one of the $k+1$ vertices on the rightmost 
path gives a size-$n$ tree with rightmost 
path of length $k$, and in this case $n$ does not end 
the rightmost path. Thus we have the recurrence
\begin{eqnarray*}
    u(n,1) & = & n(2n-3)!!  \\
    u(n,k) & = & \sum_{j=k-1}^{n-1}u(n-1,j)+(2n-k-2)u(n-1,k) \qquad n\ge k \ge 2,
\end{eqnarray*}
and the summand in 
(\ref{rightPathFromRoot}) satisfies this recurrence. The recurrence 
leads to a partial differential equation for the \gf 
$F(x,y)=\sum_{n\ge k\ge 0}u(n,k)\frac{\textrm{{\,\normalsize $x$}}^{n}}{\textrm{{\small $n$}} 
\textrm{{\footnotesize !}} } y^{k}$:
\[
(1-y)(1-2x)F_{x}(x,y)+y(1-y)F_{y}(x,y) + y^2 F(x,y) = y/\sqrt{1-2x},
\]
with solution
\[
F(x,y)= \frac{1 - (1-y)e^{y(1 - \sqrt{1-2x})} }{y\sqrt{1-2x}}.
\]
The first few values of $u(n,k)$ are
\[
\begin{array}{c|cccccc}
	n^{\textstyle{\,\backslash \,k}} & 0 & 1 & 2 & 3 & 4 & 5 \\
\hline 
        0&     1 &   &   & & & \\
	1&     0 &  1 &   & & & \\
 	2&     0 &  2 & 1 & & &   \\
	3&     0 &  9 & 5 & 1 & &  \\
	4&     0 & 60 & 35  & 9 & 1 &   \\ 
	5&     0 &525 & 315 & 90 & 14 & 1  \\
 
 \end{array}
\]
(The top left entry for $n=k=0$ in the array wants to be included in 
order to get the nice \gf.)

\textbf{Remark}\quad The summand in (\ref{rightPathFromRoot}) can be 
written somewhat more compactly by distinguishing the cases $k$ even 
or odd: the summand is 
\[
\binom{n-j-1}{j-1} \frac{(2n-2j+1)!!}{(2j+1)!!} \qquad \text{if $k=2j$ 
is even, and} 
\]
\[
\binom{n-j+1}{j} \frac{(2n-2j-1)!!}{(2j-3)!!} \qquad \text{if $k=2j-1$ is odd.}  
\]

Our interpretation of (\ref{rightPathFromRoot}) is equivalent to the 
assertion that 
the number of increasing ordered trees of size $n$ 
whose rightmost path from the root  has length $\ge k$ is 
$\binom{2n-k}{k}(2n-2k-1)!!$. We now prove this assertion bijectively. Form a ``vertex'' set 
$V=\{1,2,\ldots,n\}$ and an ``edge'' set $E=\{1,2,\ldots,n-k\}$ and 
make them disjoint by using subscripts $V$ and $E$ on 
their entries. Thus $\v\,V\cup E\,\v=2n-k$. It suffices to exhibit a bijection 
from the trees being counted to pairs consisting of a $k$-element 
subset $X$ of $V\cup E$ and a standard increasing ordered tree $T_{0}$ of size 
$n-k$, since these pairs are clearly counted by $\binom{2n-k}{k}(2n-2k-1)!!$. 
We use the tree illustrated in Fig. 1a) as a working example with $n=14$ 
and $k=5$.  
\begin{center} 

\begin{pspicture}(-8,-1)(8,5)
\psset{unit=.8cm}

\psline(-8,2)(-8,1)(-7,0)(-7,1)
\psline(-5,2)(-5,1)
\psline(-4,5)(-4,4)
\psline(-3,4)(-3,3)
\psline(-7,0)(-5,1)(-4,2)(-3,3)(-2,4)(-1,5)
\psline(-6,2)(-5,1)
\psline(-4,4)(-3,3)
\psline(-2,5)(-2,4)

\psdots(-8,2)(-8,1)(-7,0)(-7,1)(-5,2)(-5,1)(-4,5)(-4,4)(-3,4)(-4,2)(-2,4)(-1,5)(-6,2)(-3,3)(-2,5)

\rput(-7,-0.3){\textrm{{\footnotesize $0$}}}
\rput(-7,1.3){\textrm{{\footnotesize $1$}}}
\rput(-4.7,0.9){\textrm{{\footnotesize $2$}}}
\rput(-5,2.3){\textrm{{\footnotesize $3$}}}
\rput(-3.7,1.9){\textrm{{\footnotesize $4$}}}
\rput(-8.3,1){\textrm{{\footnotesize $5$}}}
\rput(-2.7,2.9){\textrm{{\footnotesize $6$}}}
\rput(-1.7,3.9){\textrm{{\footnotesize $7$}}} 
\rput(-4.3,4){\textrm{{\footnotesize $8$}}}
\rput(-2,5.3){\textrm{{\footnotesize $9$}}}
\rput(-1,5.3){\textrm{{\footnotesize $10$}}}
\rput(-3,4.3){\textrm{{\footnotesize $11$}}}
\rput(-4,5.3){\textrm{{\footnotesize $12$}}}
\rput(-8,2.3){\textrm{{\footnotesize $13$}}}
\rput(-6,2.3){\textrm{{\footnotesize $14$}}}

\rput(-4,-1.3){\textrm{{\footnotesize Fig. 1a)}}}

\psline(1,2)(1,1)(2,0)(2,1)
\psline(4,2)(4,1)
\psline(5,5)(5,4)
\psline(6,4)(6,3)
\psline[linecolor=red](2,0)(4,1)(5,2)(6,3)(7,4)(8,5)
\psline[linewidth=1.5pt,linecolor=blue](3,2)(4,1)
\psline[linewidth=1.5pt,linecolor=blue](5,4)(6,3)
\psline[linewidth=2pt,linecolor=blue](7,5)(7,4)

\psdots(1,2)(1,1)(2,0)(2,1)(3,2)(4,2)(4,1)(5,5)(5,4)(6,4)(6,3)(7,5)(7,4)

\psset{linecolor=blue}
\qdisk(5,2){3pt}
\qdisk(8,5){3pt}

\rput(2,-0.3){\textrm{{\footnotesize $0$}}}
\rput(2,1.3){\textrm{{\footnotesize $1$}}}
\rput(4.3,0.9){\textrm{{\footnotesize $2$}}}
\rput(4,2.3){\textrm{{\footnotesize $3$}}}
\rput(5.3,1.9){\textrm{{\footnotesize $4$}}}
\rput(0.7,1){\textrm{{\footnotesize $5$}}}
\rput(6.3,2.9){\textrm{{\footnotesize $6$}}}
\rput(7.3,3.9){\textrm{{\footnotesize $7$}}} 
\rput(4.7,4){\textrm{{\footnotesize $8$}}}
\rput(7,5.3){\textrm{{\footnotesize $9$}}}
\rput(8,5.4){\textrm{{\footnotesize $10$}}}
\rput(6,4.3){\textrm{{\footnotesize $11$}}}
\rput(5,5.3){\textrm{{\footnotesize $12$}}}
\rput(1,2.3){\textrm{{\footnotesize $13$}}}
\rput(3,2.3){\textrm{{\footnotesize $14$}}}

\rput(4,-1.3){\textrm{{\footnotesize Fig. 1b)}}}

\end{pspicture}
\end{center}
Call the first $k$ edges of the rightmost path from the root the 
\emph{base} edges  and their child vertices the 
\emph{base} vertices ($2,4,6,7,10$).
Say a base vertex is \emph{fertile} if it has a child that is not a base 
vertex and \emph{barren} otherwise. Highlight each barren vertex and 
the leftmost edge from each fertile vertex, and color the base edges red 
as in Fig 1b). Thus there are $k$ highlighted vertices/edges. The 
labels on the highlighted vertices are themselves the contribution from 
$V$ to the $k$-element set $X$. The red edges will be 
deleted by an iterative cut-and-paste procedure to get the required 
$(n-k)$-edge tree $T_{0}$ and then the positions of the highlighted edges in $T_{0}$
will determine the contribution from $E$ to $X$.

First, erase the barren vertices and standardize the vertex 
labeling---replace smallest by 0, next smallest by 1, and so on, to 
get the first tree in Fig. 2. 
\begin{center} 

\begin{pspicture}(-8,-1.3)(8,3.8)
\psset{unit=.7cm}

\psline(1,2)(1,1)(2,0)(4,1)(4,2)
\psline(5,4)(5,3)
\psline(6,3)(6,2)

\psline[linecolor=red](4,1)(6,2)(8,3)
\psline[linewidth=1.5pt,linecolor=blue](3,2)(4,1)
\psline[linewidth=1.5pt,linecolor=blue](5,3)(6,2)
\psline[linewidth=2pt,linecolor=blue](8,3)(8,4)

\psdots(1,2)(1,1)(2,0)(3,2)(4,2)(4,1)(5,4)(6,2)(6,3)(8,3)(8,4)

\rput(2,-0.3){\textrm{{\footnotesize $0$}}}
\rput(0.7,0.9){\textrm{{\footnotesize $3$}}}
\rput(4.3,0.9){\textrm{{\footnotesize $1$}}}
\rput(4,2.3){\textrm{{\footnotesize $2$}}}
\rput(6.3,1.9){\textrm{{\footnotesize $4$}}}
\rput(8.3,2.9){\textrm{{\footnotesize $5$}}}
\rput(4.7,2.9){\textrm{{\footnotesize $6$}}}
\rput(8,4.3){\textrm{{\footnotesize $7$}}} 
\rput(6,3.3){\textrm{{\footnotesize $8$}}}
\rput(5,4.3){\textrm{{\footnotesize $9$}}}
\rput(1,2.3){\textrm{{\footnotesize $10$}}}
\rput(3,2.3){\textrm{{\footnotesize $11$}}}


\psline(-10,2)(-10,1)(-9,0)(-7,1)(-7,2)
\psline(-9,1)(-9,0)
\psline(-6,4)(-6,3)
\psline(-5,3)(-5,2)

\psline[linecolor=red](-9,0)(-7,1)(-5,2)(-3,3)
\psline[linewidth=1.5pt,linecolor=blue](-8,2)(-7,1)
\psline[linewidth=1.5pt,linecolor=blue](-6,3)(-5,2)
\psline[linewidth=2pt,linecolor=blue](-3,3)(-3,4)

\psdots(-10,2)(-10,1)(-9,0)(-8,2)(-7,2)(-7,1)(-6,4)(-5,2)(-5,3)(-3,3)(-3,4)

\rput(-9,-0.3){\textrm{{\footnotesize $0$}}}
\rput(-9,1.3){\textrm{{\footnotesize $1$}}}
\rput(-10.3,0.9){\textrm{{\footnotesize $4$}}}
\rput(-6.7,0.9){\textrm{{\footnotesize $2$}}}
\rput(-7,2.3){\textrm{{\footnotesize $3$}}}
\rput(-4.7,1.9){\textrm{{\footnotesize $5$}}}
\rput(-2.7,2.9){\textrm{{\footnotesize $6$}}}
\rput(-6.3,2.9){\textrm{{\footnotesize $7$}}}
\rput(-3,4.3){\textrm{{\footnotesize $8$}}} 
\rput(-5,3.3){\textrm{{\footnotesize $9$}}}
\rput(-6,4.3){\textrm{{\footnotesize $10$}}}
\rput(-10,2.3){\textrm{{\footnotesize $11$}}}
\rput(-8,2.3){\textrm{{\footnotesize $12$}}}

\rput(-1,2){$\longrightarrow$}
\rput(10,2){$\longrightarrow$}

\end{pspicture}
\end{center}
\begin{center} 

\begin{pspicture}(-7,-1.8)(9,2)
\psset{unit=.7cm}

\psline(2.5,2)(4,1)(4.5,2)
\psline(3.5,2)(3.5,3)
\psline(4,1)(7,0)(9,1)(9,2)

\psline[linewidth=1.5pt,linecolor=blue](3.5,2)(4,1)
\psline[linewidth=1.5pt,linecolor=blue](4,1)(5.5,2)
\psline[linewidth=1.5pt,linecolor=blue](8,2)(9,1)

\psdots(2.5,2)(4,1)(4.5,2)(3.5,2)(3.5,3)(7,0)(9,1)(9,2)(5.5,2)(8,2)

\rput(7,-0.4){\textrm{{\footnotesize $0$}}}
\rput(3.9,0.6){\textrm{{\footnotesize $3$}}}
\rput(9.1,0.7){\textrm{{\footnotesize $1$}}}
\rput(9,2.3){\textrm{{\footnotesize $2$}}}
\rput(3.2,2.3){\textrm{{\footnotesize $4$}}}
\rput(5.6,2.3){\textrm{{\footnotesize $5$}}}
\rput(4.6,2.3){\textrm{{\footnotesize $6$}}}
\rput(3.5,3.3){\textrm{{\footnotesize $7$}}} 
\rput(2.4,2.3){\textrm{{\footnotesize $8$}}}
\rput(8,2.3){\textrm{{\footnotesize $9$}}}

\psline(-8.5,2)(-7,1)(-6.5,2)
\psline(-7.5,2)(-7.5,3)
\psline(-7,1)(-4,0)(-2,1)(-2,2)

\psline[linecolor=red](-7,1)(-5.5,2)

\psline[linewidth=1.5pt,linecolor=blue](-7.5,2)(-7,1)
\psline[linewidth=1.5pt,linecolor=blue](-3,2)(-2,1)
\psline[linewidth=2pt,linecolor=blue](-5.5,2)(-5.5,3)

\psdots(-8.5,2)(-7,1)(-6.5,2)(-7.5,2)(-7.5,3)(-4,0)(-2,1)(-2,2)(-5.5,2)(-3,2)

\rput(-4,-0.4){\textrm{{\footnotesize $0$}}}
\rput(-7.1,0.6){\textrm{{\footnotesize $3$}}}
\rput(-1.9,0.7){\textrm{{\footnotesize $1$}}}
\rput(-2,2.3){\textrm{{\footnotesize $2$}}}
\rput(-7.8,2.3){\textrm{{\footnotesize $5$}}}
\rput(-5.2,2.3){\textrm{{\footnotesize $4$}}}
\rput(-5.5,3.3){\textrm{{\footnotesize $6$}}}
\rput(-6.4,2.3){\textrm{{\footnotesize $7$}}}
\rput(-7.5,3.3){\textrm{{\footnotesize $8$}}} 
\rput(-8.6,2.3){\textrm{{\footnotesize $9$}}}
\rput(-3,2.3){\textrm{{\footnotesize $10$}}}

\rput(0.5,2){$\longrightarrow$}
\rput(0.5,-2){\textrm{{Fig. 2}}}

\end{pspicture}
\end{center}

Now eliminate the remaining red edges, one at 
a time, in increasing order of their child vertices.
Let $b$ denote the current 
smallest base vertex (initially $b=2$). Cut out the subtree rooted 
at $b$, discarding the label $b$ and its parent edge, and re-root this 
subtree at $b-1$, placed so that it lies to the right of the 
existing subtree rooted at $b-1$, and standardize vertex labels. 
Repeat this process on the second, third, \ldots, red edge in turn to obtain 
the desired $(n-k)$-edge tree $T_{0}$ (the last tree in Fig. 2) with some 
edges highlighted. Observe that the progress of individual 
highlighted edges in the successive trees of Fig. 2 can be discerned 
even though the labels on their endpoints will change.  

Finally, list the edges of $T_{0}$ in \emph{standard} order, that 
is, in increasing order of their parent vertices (preserving, of 
course, the order of edges with a common parent vertex), and take the 
positions of the highlighted edges in this list as the contribution 
from $E$ to the $k$-element set $X$. In the example, $T_{0}$ has 
defining edge list 
\[
\big( 
(0,3),(0,1),(1,9),(1,2),(3,8),(3,4),(3,6),(3,5),(4,7)\big)
\]
and the 
highlighted edges---$(1,9),(3,4),(3,5)$---are in positions 3,6,8. The 
net result is the pair $(X,T_{0})$ with 
$X=\{4_{V},10_{V},3_{E},6_{E},8_{E}\}$ and $T_{0}$ as just given.

Is this process reversible? The entries in X with subscript $E$ determine the highlighted  edges 
in $T_{0}$, and the entire process can now be reversed step-by-step 
provided we know the order in which the highlighted edges were 
originally processed. But this order is precisely their left-to-right order 
in the standard listing of the edges of $T_{0}$. \qed

\textbf{Remark}\quad $\binom{2n-k}{k}(2n-2k-1)!!$ is also the number 
of increasing ordered trees of size $n+1$ whose root has $k+1$ 
children, increasing from left to right: given an increasing ordered tree of size $n$ 
whose rightmost path from the root has length $\ge k$, delete the base 
edges, increment all labels by 1, and then attach the old root and the $k$ 
base vertices in increasing order to a new root 0. This is a bijection to the trees in question.

\subsection{}\label{maxyoung}
\vspace*{-7mm}
\begin{equation}
 \sum_{k=1}^{n}\frac{(n-1)!}{(k-1)!}\,k\,(2k-3)!! = (2n-1)!!
    \label{youngleaf}
\end{equation}
This identity counts increasing ordered trees of size $n$ by the maximum 
$k$ of the
young leaves where a \emph{young} leaf is a leaf with no left sibling (leaf 
or otherwise). Every nonempty  increasing ordered tree has 
at least one young leaf. To establish this count, consider the $2n-1$ 
ways to produce an increasing ordered tree of $n$ edges by inserting 
a leaf $n$ into an increasing ordered tree of $n-1$ 
edges: either as the leftmost child of one of the $n$ vertices or so 
that the new edge is 
the immediate right neighbor edge of one of the $n-1$ existing edges.
In the first case, the maximum young leaf becomes $n$; in the second, 
it is preserved. In particular, the maximum young leaf becomes $n$ 
in $n$ ways from \emph{any} increasing ordered tree of $n-1$ 
edges. Hence, if $u(n,k)$ denotes the number of increasing ordered trees 
with $n$ edges whose maximum young leaf is $k$, we have 
$u(n,n)=n(2n-3)!!$ and the recurrence $u(n,k)=(n-1)u(n-1,k)$ for $1\le 
k <n$. Iterating the recurrence yields 
$u(n,k)=(n-1)^{\underline{n-k}}\,k\,(2k-3)!!$ for all $1\le k \le n$.

The bivariate \gf $\sum_{n,k\ge 1}\frac{\textrm{{\small $(n-1)$}} 
\textrm{{\footnotesize !}} }{\textrm{{\small $(k-1)$}} 
\textrm{{\footnotesize !}} }\,k\,(2k-3)!!
\frac{\textrm{{\,\normalsize $x$}}^{n-1}}{\textrm{{\small $(n-1)$}} 
\textrm{{\footnotesize !}} }
y^{k-1}$ is 
\[
 \frac{1-xy}{(1-x)(1-2xy)^{3/2}},
\]
and the first few values are
\[
\begin{array}{c|ccccc}
	n^{\textstyle{\,\backslash \,k}} & 1 & 2 & 3 & 4 & 5 \\
\hline 
	1&    1 &   & & & \\
 	2&    1 & 2 & & &   \\
	3&    2 & 4 & 9 & &  \\
	4&    6 & 12  & 27 & 60 &   \\ 
	5&    24 & 48 & 108 & 240 & 525  \\
 
 \end{array}
\]

\subsection{}
\vspace*{-7mm}
\begin{equation}
   \sum_{k=0}^{n/2}\binom{n}{2k}\binom{2k}{k}\frac{n!}{2^{2k}}=(2n-1)!!
    \label{pm1}
\end{equation}
This identity counts perfect matchings of $[2n]$ by number $k$ of 
matches in which both entries are $\leq n$. In fact, (\ref{pm1}) is 
the special case $r=0$ of a family of identities indexed by a nonnegative 
integer $r$:
\begin{equation}
\sum_{k=r}^{(n+r)/2}\binom{n}{2k-r}
   \binom{2k-r}{k}(n+r)^{\underline{r}}\,(n-r)!\frac{1}{2^{2k-r}}=(2n-1)!!
    \label{pm2}
\end{equation}
where $(n+r)^{\underline{r}}=(n+r)(n+r-1)\ldots(n+1)$ is the falling 
factorial. For given $r$, as a straightforward direct count shows, 
(\ref{pm2}) counts perfect matchings of $[2n]$ by number $k$ of 
matches in which both entries are $\leq n+r$ (also by number $k$ of 
matches in which both entries are $\leq n-r$). 
The bivariate \gf $\sum_{n,k\ge r}\binom{2k-r}{k}(n+r)^{\underline{r}}\,
(n-r)!\frac{1}{\textrm{{\footnotesize \raisebox{-1mm}{2}}} ^{2k-r}}\:
\frac{\textrm{{\,\normalsize $x$}}^{n-r}}{\textrm{{\small $(n-r)$}} 
\textrm{{\footnotesize !}} }\:
y^{k-r}$ is 
\[
\frac{(2r-1)!!}{\big( (1-x)^2-x^2 y \big)^{ (2r+1)/2 }}.
\]
For $r=0$, the first few values  are
\[
\begin{array}{c|ccccc}
	n^{\textstyle{\,\backslash \,k}} & 0 & 1 & 2 & 3   \\
\hline 
	0&    1 &   & &  \\
        1&    1 &   & &  \\
 	2&    2 & 1& &    \\
	3&    6 & 9  & &   \\
	4&    24 & 72  & 9 &    \\ 
	5&    120 & 600 & 225 &    \\
	6&    720 & 5400 & 4050 & 225    \\

 \end{array}
\]
and for $r=1$, the first few values  are
\[
\begin{array}{c|cccc}
	n^{\textstyle{\,\backslash \,k}}  & 1 & 2 & 3   \\
\hline 
        1&    1 &   &   \\
 	2&    3 & &     \\
	3&    12 & 3  &    \\
	4&    60 & 45  &      \\ 
	5&    360 & 540 & 45     \\
	6&    2520 & 6300 & 1575    \\

 \end{array}
\]

By transposing factors, 
(\ref{pm2}) can be written in the alternative form
\begin{equation}
\sum_{k=r}^{(n+r)/2}\binom{n}{2k-r}
   \binom{2k-r}{k}2^{n-2k+r}=\binom{2n}{n-r},
    \label{pm3}
\end{equation}
an identity that counts bicolored \emph{UDF} paths of length $n$ ending at height $r$ by 
number $k$ of upsteps. A  \emph{UDF} path is a lattice path of upsteps $(1,1)$, downsteps $(1,-1)$, and 
flatsteps $(1,0)$; bicolored means each flatstep is colored red or 
blue; height is relative to the horizontal line through the initial 
vertex. 

The case $r=1$ of (\ref{pm2}) also counts 
increasing ordered trees of $n$ edges by number $k$ of 
young leaves (as defined in Section \ref{maxyoung}). To see this, note that 
inserting a leaf $n$ into an increasing ordered tree of $n-1$ edges 
($2n-1$ possible ways) always either preserves or increments (by 1) the number 
of young leaves, and the number of ways to increment is 
$n-2k$ where $k$ is the number of young leaves. This observation leads to the recurrence
\[
u(n,k)=(n+2k-1)u(n-1,k) + (n-2k+2)u(n-1,k-1) \qquad  1\le k \le 
(n+1)/2
\]
for the number $u(n,k)$ of increasing ordered trees of $n$ edges 
and $k$ young leaves, and the recurrence is satisfied by the summand.

\subsection{}\label{stirliByFirst}
\vspace*{-7mm}
\begin{equation}
\sum_{k=1}^{n}\frac{ (2n-2)!! (2k-3)!! }{(2k-2)!!} = (2n-1)!!
    \label{StirlingByFirst}
\end{equation}
This identity counts 

\noindent(1) Stirling permutations of size $n$ by first entry $k$, 

\noindent(2) Stirling permutations of size $n$ by position $2k-1$ of the first 1 (the
position of the first 1 is necessarily odd),

\noindent(3) increasing ordered trees of size $n$ by the parent $k-1$ of $n$, 

\noindent(4) increasing ordered trees of size $n$ by the label $k$ on the leaf in 
the minimal path from the root. The \emph{minimal path} starts at 
the root and successively travels to the smallest child vertex until 
it arrives at a leaf. 

\noindent \textbf{Proofs} 

\noindent(1) and (4)
The number of choices in building up the 
object, inserting a pair $ii$ in the permutation or a 
leaf $i$ in the tree for $i=1$ to $n$, is successively 
$1,3,\ldots,2k-3,1,2k,2k+2,\ldots,2n-2$ and their product is the 
summand.

\noindent(2) If a letter occurs to the
left of the first 1, then both occurrences do so, and so the number of such
permutations is $\binom{n-1}{k-1}$ [choose support set for the first $2k-2$
entries] $\times (2k-3)!!$ [form a Stirling permutation on this
support] $\times (2n-2k)!!$ [form a Stirling permutation of size $n-k+1$
that starts with a 1], and $\binom{n-1}{k-1}(2k-3)!!(2n-2k)!!$ is
an equivalent expression for the summand.

\noindent(3) Here, build up the tree by successively inserting leaves 
$1,2,\ldots,k-1$, then $n$ must be inserted as a child of $k-1$, then 
proceed to insert $k,k+1,\ldots,n-1$. The number of choices is the 
same as in (1) and (4) above.

The bivariate \gf $\sum_{n\ge 1,k\ge 1} 
\frac{\textrm{{\footnotesize $(2n-2)!! (2k-3)!!$}}}{\textrm{{\footnotesize 
$(2k-2)!!$}}}\, 
\frac{\textrm{{\,\normalsize $x$}}^{n-1}}{\textrm{{\small $(n-1)$}} 
\textrm{{\footnotesize !}} } y^{k-1}$  is
\[
\frac{1}{(1-2x)\sqrt{1-2x y}},
\]
 and the first few values are
\[
\begin{array}{c|ccccc}
	n^{\textstyle{\,\backslash \,k}} & 1 & 2 & 3 & 4 & 5\\
\hline 
	1&    1 &   & & & \\
 	2&    2 & 1 & & &   \\
	3&    8 & 4  & 3 & &  \\
	4&    48 & 24  & 18 & 15 &   \\ 
	5&    384 & 192 & 144 & 120 & 105  \\

 \end{array}
\]
With rows reversed, this array is entry
\htmladdnormallink{A122774}{http://www.research.att.com:80/cgi-bin/access.cgi/as/njas/sequences/eisA.cgi?Anum=A122774} in
OEIS \cite{oeis}, and the reversed array counts (i) height-labeled Dyck paths by the 
position among the upsteps of the last upstep with label 1, and (ii) increasing ordered trees 
by the maximum child of $1$.

\subsection{}
\vspace*{-7mm}
\begin{equation}
(2n-2)!! + \sum_{k=2}^{n} \frac{(2n-1)!!(2k-4)!!}{(2k-1)!!} = (2n-1)!!
    \label{MinChildOf1}
\end{equation}
This identity counts \IOTs of size $n$ by smallest
child $k$ of 1 ($k=1$ if vertex 1 has no children). To see 
this, let 
$u(n,k)$ be the number of such trees. Consideration of the effect of inserting 
a leaf $n$ into an \IOT of size $n-1$ on the smallest
child of 1 yields the recurrence
\begin{eqnarray*}
  u(n,1) & = & (2n-2)!!  \\
  u(n,k) & = & (2n-1)u(n-1,k) \qquad  2\le k \le n-1  \\
  u(n,n) & = & (2n-4)!!,
\end{eqnarray*}
satisfied by the summands. The recurrence
leads to a first-order ordinary differential equation for the \gf 
$F(x,y)=\sum_{n \ge k \ge 1}u(n,k)$ $\frac{\textrm{{\,\normalsize $x$}}^{n-1}}{\textrm{{\small $(n-1)$}} 
\textrm{{\footnotesize !}} } y^{k-1}$:
\[
(1-2x)F_{x} - 3F + 1/(1-2x) = \frac{y}{1 - 2xy},
\]
with solution
\[
F(x,y)=\frac{1}{(1-2x)^{3/2}}+\frac{\sqrt{y-y^{2}}}{(1-2x)^{3/2}y}\,
\tan^{-1}\left(\frac{\sqrt{y-y^{2}}\,(-1+2xy+\sqrt{1-2x})}{1-2y+2xy^{2}}\right).
\]
The first few values of $u(n,k)$ are
\[
\begin{array}{c|ccccc}
	n^{\textstyle{\,\backslash \,k}} & 1 & 2 & 3 & 4 & 5 \\
\hline 
	1&    1 & & & &\\
 	2&    2 & 1 & & &  \\
	3&    8 & 5 & 2 &  \\
	4&    48  & 35 & 14 & 8 &  \\ 
	5&    384  & 315 & 126 & 72 & 48 \\

 \end{array}
\]

With rows reversed, this array has an interpretation for which the recurrence 
relation is not immediately obvious: let $v(n,k)$ denote the number of
Stirling permutations of size $n$ for which the maximum $M$
of the entries preceding the first 1 (taken as $n+1$ if the 
permutation starts with a 1) is $k$  ($2\le k \le n+1$). Then $v(n,k)=u(n,n+2-k)$. 
This follows from the next Proposition by comparing recurrences.
\begin{prop}\label{maxRec}
\begin{eqnarray*}
      (i)  & \ v(n,n+1) & =\quad(2n-2)!!, \\
     (ii)  & \ v(n,2)   & =\quad(2n-4)!!,\textrm{\quad and}  \\
    (iii)  & \ v(n,k)   & =\quad(2n-1)v(n-1,k-1) \textrm{\quad for $3\le k \le n$.}
\end{eqnarray*}
 \end{prop}
\noindent \textbf{Proof}

(i) \ Consider permutations that start with a 1. Inserting two 
adjacent $n$'s immediately after any one of the $2n-2$ entries in such a permutation of 
size $n-1$ gives one of size $n$. Thus $v(n,n+1)=(2n-2)v(n-1,n)$, 
and (i) follows.

(ii) \  A Stirling permutation with $M=2$ necessarily begins 
221\ldots. Deleting the initial 2s is a bijection to Stirling 
permutations of size $n-1$ that start with a 1, counted by $v(n-1,n)=(2n-4)!!$.

(iii) \  Now suppose $3 \le k \le n.$ A Stirling permutation of size $n-1$ with $M=k-1$ 
yields $2n-1$ Stirling permutations of size $n$ with $M=k$ (all distinct) as 
follows. Increment each entry $\ge 2$ by 1  and tentatively insert a 
pair of adjacent 2s anywhere in the resulting permutation ($2n-1$ 
choices) to obtain a permutation $\si$. Let $I$ denote the initial 
segment of $\si$ terminating at the first 2. If no entry exceeding 2 occurs 
exactly once in $I$, then $\si$ is already a Stirling permutation and 
we may leave the tentative 2s in place. Otherwise, choose the 
smallest $i>2$ that occurs exactly once in $I$, say $i_{1}$, and 
interchange the 2s and $i_{1}$s in $\si$ to obtain a perm $\si_{1}$.
Let $I_{1}$ denote the initial 
segment of $\si_{1}$ terminating at the first $i_{1}$. If no entry 
exceeding $i_{1}$ occurs 
exactly once in $I_{1}$, then $\si_{1}$ is a Stirling permutation, and 
stop. Otherwise, proceed similarly to get $i_{2}>i_{1}$, interchange the 
$i_{1}$s and $i_{2}$s to get $\si_{2}$ and continue until you arrive 
at a permutation $\si_{k}$ that does have the Stirling property. The 
original Stirling permutation of size $n-1$ and the location of the 
tentative 2s can be recovered from $\si_{k}$ and so this process is a 
bijection $\V(n-1,k-1)\times [2n-1] \rightarrow \V(n,k)$ where $\V(n,k)$ 
is the set counted by $v(n,k)$. \qed

Also, $v(n,k)$ is the number of \IOTs of size $n$ with $k$ the maximal 
descendant of 1 (taken as $n+1$ if 1 is a leaf since ``descendant'' 
here means ``proper descendant''). The \gf $\sum_{n\ge 1}\sum_{k=2}^{n+1}v(n,k) 
\frac{\textrm{{\,\normalsize $x$}}^{n-1}}{\textrm{{\small $(n-1)$}} 
\textrm{{\footnotesize !}} } y^{k-2}$ for $v(n,k)$ is marginally more 
concise than 
that for $u(n,k)$:
\[
\frac{1}{(1-2xy)^{3/2}}+\frac{\sqrt{y-1}}{(1-2xy)^{3/2}}\,
\tan^{-1}\left(\frac{\sqrt{y-1}\,(1-2x-\sqrt{1-2xy})}{2-2x-y}\right).
\]

\subsection{}
\vspace*{-7mm}
\begin{equation}
(2n-3)!! + \sum_{k=1}^{n-1}\frac{2\,(2n-1)!!}{(2k+1)(2k-1)} = (2n-1)!!
    \label{minAftern}
\end{equation}
This identity counts Stirling permutations of size $n$ by smallest 
entry $k$ following the last $n$ ($k=0$ if the last entry is $n$). The 
recurrence relation for these permutations is
\begin{eqnarray*}
  u(n,0) & = & (2n-3)!!  \\
  u(n,k) & = & (2n-1)u(n-1,k) \qquad  1\le k \le n-2  \\
  u(n,n-1) & = & 2(2n-5)!!,
\end{eqnarray*}
satisfied by the summand.

\noindent \textbf{Remark}\quad The identity itself is trivial to prove, since 
the sum is telescoping, but the \gf 
$\sum_{n \ge 1,\: 0 \le k \le n-1}u(n,k)\frac{\textrm{{\,\normalsize $x$}}^{n-1}}{\textrm{{\small $(n-1)$}}
\textrm{{\footnotesize !}} } y^{k}$ is cumbersome:
\begin{multline*}
\frac{1}{(1 - 2x)^{3/2}} + \frac{1}{\sqrt{1 - 2x}} - \frac{\sqrt{1 - 
2xy}}{1 - 2x} - \\
 \frac{1 - y}{2\sqrt{y}  (1 - 2x)^{3/2}} 
\log\left( \frac{1 + y - 4xy - 2\sqrt{(1 - 2x)y(1 - 2xy)}}{(1 - 
\sqrt{y})^2}\right). 
\end{multline*} 
The first few values of $u(n,k)$ are
\[
\begin{array}{c|ccccc}
	n^{\textstyle{\,\backslash \,k}} & 0 & 1 & 2 & 3 & 4  \\
\hline 
	1&    1 &   & & & \\
 	2&    1 & 2 & & &   \\
	3&    3 & 10  & 2 & &  \\
	4&    15 & 70  & 14 & 6 &   \\ 
	5&    105 & 630 & 126 & 54 & 30  \\
\end{array}
\]

\subsection{}\label{firstAscentLength}
\vspace*{-7mm}
\begin{equation}
\sum_{k=1}^{n}k! \binom{2n-k-1}{k-1}\,(2n-2k-1)!! = (2n-1)!!
    \label{firstascent}
\end{equation}

This identity counts increasing ordered trees of $n$ edges by
outdegree $k$ of the root. The $k!$ factor allows us to assume that the 
$k$ children of the root increase from left to right, and the number 
of such trees is given by the Remark at the end of Section \ref{rightPath}.

The identity also counts height-labeled Dyck $n$-paths by length $k$ of first 
ascent. To show this, we exhibit a bijection from height-labeled Dyck 
$n$-paths to increasing ordered $n$-trees that sends ``length of first 
ascent'' to ``number of children of the root''. The tree is built up 
using trees with some of their leaves unlabeled. The construction involves 
a pair of 
sequences $\big(a(i)\big)_{i=1}^{n}$ and $\big(b(i)\big)_{i=1}^{n}$ 
that characterizes the path:
\begin{center} \vspace*{-2mm}
\begin{tabular}{l}
  $a(i)= \#\,$upsteps immediately preceding the $i$th downstep, and  
  \\[1mm]
  $b(i)= $ label on the upstep matching the $i$th downstep.  \\
  \end{tabular}    
\end{center}
Thus $0\le a(i) \le n$ and $\sum_{i=1}^{n}a(i)=n$. At step $i$ ($1\le i 
\le n$), attach $a(i)$ unlabeled 
leaves to (leaf) vertex $i-1$ and apply label $i$ to the $b(i)$-th 
unlabeled leaf in walkaround order as illustrated. 
\Einheit=0.6cm
\[
\Pfad(-13,0),33433344344434\endPfad
\gray{\Pfad(-13,0),11111111111111\endPfad }
\DuennPunkt(-13,0)
\DuennPunkt(-12,1)
\DuennPunkt(-11,2)
\DuennPunkt(-10,1)
\DuennPunkt(-9,2)
\DuennPunkt(-8,3)
\DuennPunkt(-7,4)
\DuennPunkt(-6,3)
\DuennPunkt(-5,2)
\DuennPunkt(-4,3)
\DuennPunkt(-3,2)
\DuennPunkt(-2,1)
\DuennPunkt(-1,0)
\DuennPunkt(0,1)
\DuennPunkt(1,0)
\Label\o{\textrm{{\footnotesize $\xleftarrow{\text{\phantom{aaa}}}$ 
matching steps
$\xrightarrow{\text{\phantom{aaa}}}$ }}}(-5.9,1)
\Label\o{\textrm{{\scriptsize 1}}}(-12.7,0.3)
\Label\o{\textrm{{\scriptsize 2}}}(-11.7,1.3)
\Label\o{\textrm{{\scriptsize 1}}}(-9.7,1.3)
\Label\o{\textrm{{\scriptsize 3}}}(-8.7,2.3)
\Label\o{\textrm{{\scriptsize 3}}}(-7.7,3.3)
\Label\o{\textrm{{\scriptsize 2}}}(-4.7,2.3)
\Label\o{\textrm{{\scriptsize 2}}}(-11.7,1.3)
\Label\o{\textrm{{\scriptsize 1}}}(-0.7,0.3)
\Label\o{ \textrm{$ 
\begin{array}{c|ccccccc}
    i & 1 & 2 & 3 & 4 & 5 & 6 & 7  \\ \hline
    a(i) & 2 & 3 & 0 & 1 & 0 & 0 & 1  \\
    b(i) & 2 & 3 & 3 & 2 & 1 & 1 & 1
\end{array} $}}(8,1.2)
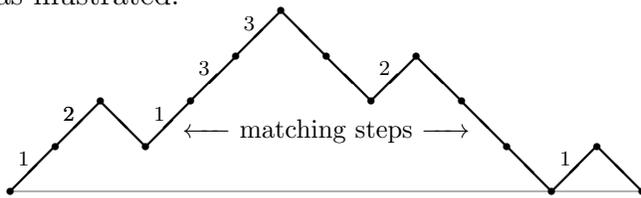
\Label\o{\textrm{{\small height-labeled Dyck 7-path}}}(-6,-1.5)
\Label\o{\textrm{{\small defining sequences}}}(8.3,-1.5)
\]

\begin{center} 
\begin{pspicture}(-7.2,-1.8)(8.8,2)
\psset{unit=.8cm}

\psline(-8,1)(-7,0)(-6,1)
\psline(-4,1)(-3,0)(-2,1)(-3,2)
\psline(-2,2)(-2,1)(-1,2)
\psline(1,1)(2,0)(3,1)(2,2)
\psline(3,2)(3,1)(4,2)
\psline(6,1)(7,0)(8,1)(7,2)
\psline(8,3)(8,2)(8,1)(9,2)

\psdots(-8,1)(-7,0)(-6,1)(-4,1)(-3,0)(-2,1)(-3,2)(-2,2)(-1,2)(1,1)(2,0)(3,1)(2,2)(3,2)(4,2)(6,1)(7,0)(8,1)(7,2)(8,3)(8,2)(9,2)

\rput(-6,1.3){\textrm{{\footnotesize $1$}}}
\rput(-1.7,0.8){\textrm{{\footnotesize $1$}}}
\rput(3.3,0.8){\textrm{{\footnotesize $1$}}}
\rput(8.3,0.8){\textrm{{\footnotesize $1$}}}
\rput(-2,2.3){\textrm{{\footnotesize $2$}}}
\rput(4,2.3){\textrm{{\footnotesize $3$}}}
\rput(9,2.3){\textrm{{\footnotesize $3$}}}
\rput(3,2.3){\textrm{{\footnotesize $2$}}}
\rput(8.3,2.3){\textrm{{\footnotesize $2$}}}
\rput(7,2.3){\textrm{{\footnotesize $4$}}}

\rput(-7,-0.3){\textrm{{\footnotesize $0$}}}
\rput(-3,-0.3){\textrm{{\footnotesize $0$}}}
\rput(2,-0.3){\textrm{{\footnotesize $0$}}}
\rput(7,-0.3){\textrm{{\footnotesize $0$}}}

\rput(-7,-1.3){\textrm{{\small step 1}}}
\rput(-3,-1.3){\textrm{{\small step 2}}}
\rput(2,-1.3){\textrm{{\small step 3}}}
\rput(7,-1.3){\textrm{{\small step 4}}}

\rput(-5,1.5){$\rightarrow$}
\rput(0,1.5){$\rightarrow$}
\rput(5,1.5){$\rightarrow$}

\end{pspicture}
\end{center}
\begin{center} 
\begin{pspicture}(-7,-2.5)(7,3)
\psset{unit=.8cm}

\psline(-7,1)(-6,0)(-5,1)(-6,2)
\psline(-5,3)(-5,2)(-5,1)(-4,2)
\psline(-1,1)(0,0)(1,1)(0,2)
\psline(1,3)(1,2)(1,1)(2,2)
\psline(5,1)(6,0)(7,1)(6,2)
\psline(7,4)(7,3)(7,2)(7,1)(8,2)

\psdots(-7,1)(-6,0)(-5,1)(-6,2)(-5,3)(-5,2)(-4,2)(-1,1)(0,0)(1,1)(0,2)(1,2)(1,3)(2,2)
(5,1)(6,0)(7,1)(6,2)(7,3)(7,4)(7,2)(8,2)

\rput(-4.7,0.8){\textrm{{\footnotesize $1$}}}
\rput(1.3,0.8){\textrm{{\footnotesize $1$}}}
\rput(7.3,0.8){\textrm{{\footnotesize $1$}}}
\rput(-4,2.3){\textrm{{\footnotesize $3$}}}
\rput(2,2.3){\textrm{{\footnotesize $3$}}}
\rput(8,2.3){\textrm{{\footnotesize $3$}}}
\rput(-4.7,2.3){\textrm{{\footnotesize $2$}}}
\rput(1.3,2.3){\textrm{{\footnotesize $2$}}}
\rput(7.3,2.3){\textrm{{\footnotesize $2$}}}
\rput(6,2.3){\textrm{{\footnotesize $4$}}}
\rput(0,2.3){\textrm{{\footnotesize $4$}}}
\rput(-6,2.3){\textrm{{\footnotesize $4$}}}
\rput(-7,1.3){\textrm{{\footnotesize $5$}}}
\rput(-1,1.3){\textrm{{\footnotesize $5$}}}
\rput(5,1.3){\textrm{{\footnotesize $5$}}}
\rput(1,3.3){\textrm{{\footnotesize $6$}}}
\rput(7.3,3.3){\textrm{{\footnotesize $6$}}}
\rput(7,4.3){\textrm{{\footnotesize $7$}}}

\rput(-6,-0.3){\textrm{{\footnotesize $0$}}}
\rput(0,-0.3){\textrm{{\footnotesize $0$}}}
\rput(6,-0.3){\textrm{{\footnotesize $0$}}}

\rput(-6,-1.3){\textrm{{\small step 5}}}
\rput(0,-1.3){\textrm{{\small step 6}}}
\rput(6,-1.3){\textrm{{\small step 7}}}

\rput(-8.5,1.5){$\rightarrow$}
\rput(-2.5,1.5){$\rightarrow$}
\rput(3.5,1.5){$\rightarrow$}

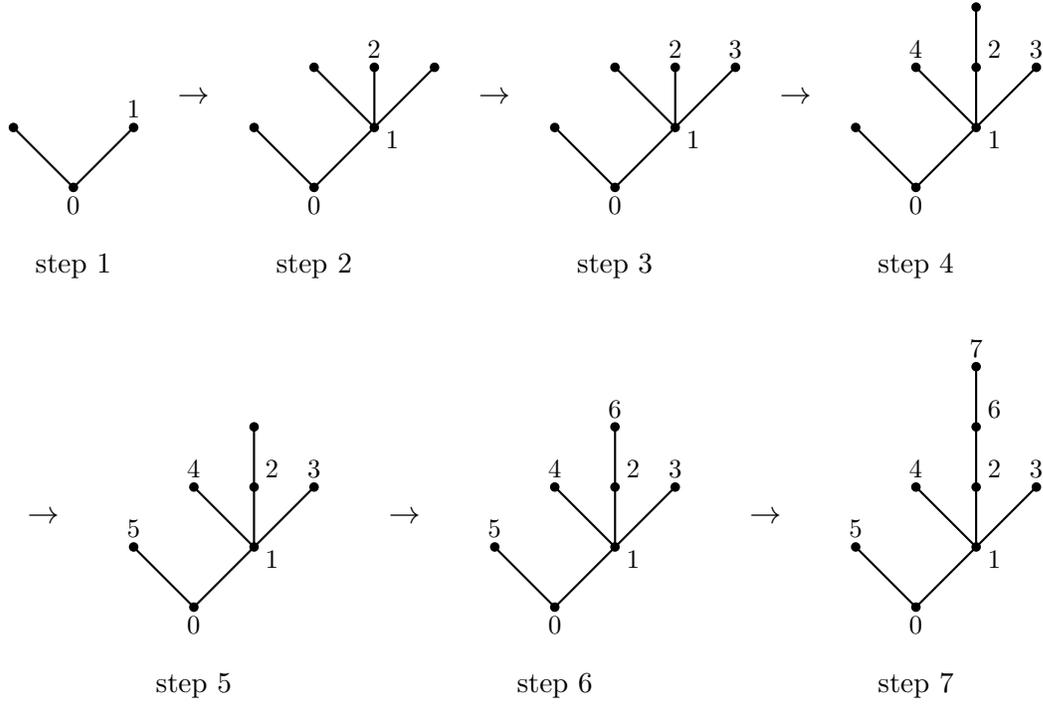
\rput(0,-2.6){\textrm{{\small construction of corresponding tree}}}

\end{pspicture}
\end{center}

We remark without proof that the identity also counts leaf-labeled 0-2 trees of size $n$ by length 
of the path from $1$ to the root.

The bivariate \gf $\sum_{n\ge 0,k\ge 0} k! \binom{2n-k-1}{k-1}$ $(2n-2k-1)!!  
\frac{\textrm{{\normalsize $x$}}^{n}}{\textrm{{\small $n$}} 
\textrm{{\footnotesize !}} } y^{k}$ is 
\[
\frac{1-y-y\sqrt{1-2x}}{1-2y+2xy^{2}}.
\]
Omitting the empty tree ($n=k=0$),  the \gf $\sum_{n\ge 1,k\ge 1} k! \binom{2n-k-1}{k-1}$ $(2n-2k-1)!!  
\frac{\textrm{{\normalsize $x$}}^{n-1}}{\textrm{{\small $(n-1)$}} 
\textrm{{\footnotesize !}} } y^{k-1}$ looks rather different: 
\[
\frac{1}{\sqrt{1-2x}\:\big(1-y+y\sqrt{1-2x}\,\big)^{2}},
\]
and the first few values are
\[
\begin{array}{c|ccccc}
	n^{\textstyle{\,\backslash \,k}} & 1 & 2 & 3 & 4 & 5 \\
\hline 
	1&    1 &   & & & \\
 	2&    1 & 2 & & &   \\
	3&    3 & 6  & 6 & &  \\
	4&    15 & 30  & 36 & 24 &   \\ 
	5&    105 & 210 & 270 & 240 & 120  \\

\end{array}
\]
This array is entry
\htmladdnormallink{A102625}{http://www.research.att.com:80/cgi-bin/access.cgi/as/njas/sequences/eisA.cgi?Anum=A102625}
in OEIS.

The result of Section \ref{firstAscentLength} can be refined to give 
the joint  distribution of first ascent 
length and first peak upstep label: let $u(n,j,k)$ denote the number 
of height-labeled Dyck $n$-paths 
whose first peak upstep has label $j$ and whose
initial ascent has length $k\ (n\ge k\ge j\ge 1)$.
Then, since $j$ is uniformly distributed over $[k]$,
$u(n,j,k)=(k-1)! \binom{2n-k-1}{k-1}\,(2n-2k-1)!!$ for $1\le j \le k$,
leading to the \gf $\sum_{n\ge 1,j\ge 1,k\ge 1}u(n,j,k)\frac{\textrm{{\normalsize $x$}}^{n-1}}{\textrm{{\small $(n-1)$}} 
\textrm{{\footnotesize !}} }y^{j-1}z^{k-1}=$
\[
\frac{1}{\sqrt{1-2x}\,\big(1-yz+yz\sqrt{1-2x}\,\big)\:\big(1-z+z\sqrt{1-2x}\,\big)}.
\]
The bijection of this section shows that $u(n,j,k)$ also counts \IOTs whose 
root has $k$ children among which vertex 1 is the $j$-th.

\section{Non-Round Identities}\label{other}
\vspace*{-5mm}
\subsection{}\label{firstdescentsubsection}
\vspace*{-7mm}
\begin{equation}
n!+ \sum_{k=1}^{n-1}\big((k-1)!! (2n-k)!! -  k!!  (2n-k-1)!!\big)= (2n-1)!!
    \label{firstdescent}
\end{equation} 
This identity is trivial to prove---the sum is telescoping and collapses to
$(2n-1)!!-$ \linebreak
$(n-1)!!n!!$---but it has an interesting 
interpretation: it counts height-labeled Dyck $n$-paths by length $k$ 
of the first 
descent where a \emph{descent} is a maximal sequence of contiguous downsteps 
and the term $n!$ corresponds to $k=n$. To see this, let 
$u(n,k)$ denote the number of height-labeled Dyck $n$-paths whose first descent 
has length $\ge k$ ($1\le k \le n$). The next-size-up construction 
described in 
(\ref{htdyck}) yields one of these paths precisely when the 
$(n-1)$-path has first descent of length $\ge k$ and the specified 
vertex is not the terminal vertex of one of the first $k-1$ downsteps 
in the first descent. Thus $u(n,k)=\big( 
(2n-1)-(k-1)\big)u(n-1,k)=(2n-k)u(n-k)$ for $n\ge k$. Together with 
the obvious initial case $u(k,k)=k!$, this recurrence yields that 
$u(n,k)=k!(k+2)^{\overline{n-k,2}}$, where 
$k^{\overline{n,2}}=k(k+2)(k+4)\cdots$ to $n$ factors is the rising 
double factorial. Equivalently, $u(n,k)=(k-1)!!(2n-k)!!$. Thus the 
number of height-labeled Dyck $n$-paths with first 
descent of length $k$ is $u(n,n)=n!$ for $k=n$, and  
$u(n,k)-u(n,k+1) = (k-1)!! (2n-k)!! -  k!!  (2n-k-1)!!$ for $1 \le 
k\le n-1$, and the first few values are
\[
\begin{array}{c|ccccc}
	n^{\textstyle{\,\backslash \,k}} & 1 & 2 & 3 & 4 & 5 \\
\hline 
	1&    1 &   & & & \\
 	2&    1 & 2 & & &   \\
	3&    7 & 2  & 6 & &  \\
	4&    57 & 18  & 6 & 24 &   \\ 
	5&    561 & 174 & 66 & 24 & 120  \\

 \end{array}
\]

\subsection{}\label{gessel}
\vspace*{-7mm}
\begin{equation}
\sum_{k=1}^{n} \SecondEulerian{n}{k} = (2n-1)!!
    \label{eulerian2}
\end{equation}
Here $\SecondEulerianInline{n}{k}$ is 
the second-order Eulerian number (indexed so that $1\le k \le n$)
\htmladdnormallink{A008517}{http://www.research.att.com:80/cgi-bin/access.cgi/as/njas/sequences/eisA.cgi?Anum=A008517}.
This identity counts 

\noindent(1) Stirling permutations of size $n$ by number $k$ of 
descents (including a conventional descent at the end), 

\noindent(2) Stirling permutations of size $n$ by number $k$ of plateaus, 
that is, pairs of adjacent equal entries,

\noindent(3) increasing ordered trees of size $n$ by number $k$ of 
leaves, 

\noindent(4) height-labeled Dyck $n$-paths by number $k$ of upstep-free vertices, and

\noindent(5) trapezoidal words of length $n$ by number $k$ of distinct entries. 

\noindent Furthermore, the second-order Eulerian triangle \emph{with 
reversed rows} counts 

\noindent(6) \IOTs by number of descents where a descent in a tree is a 
pair of adjacent sibling vertices with the first larger than the 
second (no conventional descents), and

\noindent(7) height-labeled Dyck $n$-paths by number of peaks.

\noindent \textbf{Proofs} 

\noindent(1) Several proofs are given in \cite{stirpoly}.

\noindent(2) Stirling permutations are usually defined on support set $[n]$ 
but of course can be similarly defined on an arbitrary set of positive 
integers. We present a bijection $\phi$, actually an involution, on 
Stirling permutations of arbitrary support set that preserves size and 
interchanges  ``\#\:descents'' and ``\#\:plateaus''. First, $\phi$ is 
the identity on the empty permutation. A nonempty Stirling 
permutation can be written as the concatenation $A\,m\,B\,m\,C$ 
where $m$ is the smallest integer in its support set and $A,B,C$ are 
perforce themselves Stirling permutations. Now define $\phi$ 
recursively by
\[
\phi(A\,m\,B\,m\,C) = \phi(A)\,m\,\phi(C)\,m\,\phi(B).
\]
One checks, using induction, that $\phi$ has the required properties.

\noindent(3) This result follows from Janson's bijection (see Section 
\ref{incrTree}) since 
it clearly sends leaves to plateaus.

\noindent(4) There are $2n+1$ vertices in a Dyck $n$-path and a vertex is 
\emph{upstep-free} if it is not incident with an upstep. 
The recurrence of 
(\ref{htdyck}) for height-labeled Dyck paths can easily be refined to a  
recurrence for $h(n,k)$, the number of height-labeled Dyck $n$-paths with 
$k$ upstep-free vertices, which 
turns out to be $h(n,k)=k\,h(n-1,k)+(2n-k)h(n-1,k-1)$---the defining 
recurrence for second-order Eulerian numbers \cite[p.\:27]{stirpoly}.

The recurrence leads to a 
bijection from height-labeled Dyck $n$-paths to Stirling 
permutations that sends \#\,upstep-free vertices to \#\,conventional 
descents. The bijection uses identical ``codings'' for the two classes 
of objects. 
Consider the alphabet consisting of two copies of the positive 
integers, distinguished by subscripts $Y$ and $N$. Let $\C_{n}$ 
denote the set of words $\big(w(i)\big)_{i=1}^{n}$ over this alphabet 
satisfying the following condition for each $i$. 
If $w(i)$ has subscript $Y$, then its value is $\le 1+\#\,$subscripts $N$ preceding
$w(i)$, while if $w(i)$ has subscript $N$, then its value is $\le 2i-2-\#\,$subscripts $N$ preceding
$w(i)$.
Clearly, there are $2i-1$ choices for $w(i)$ regardless of the 
preceding entries, and so $\v\,\C_{n}\,\v=(2n-1)!!$. Indeed, a word 
$w$ in $\C_{n}$ can be represented as a lattice path along with a set 
of lattice points: discard $w(1)$---necessarily $1_{Y}$---and subtract 1 
from each remaining value to get a sequence of $n-1$ nonnegative 
integers $\big(b(i)\big)$ and a sequence of $n-1$ subscripts. The subscripts 
give the path letting $Y$ 
denote a downstep $(1,-1)$ and $N$ a flatstep $(1,0)$, and the 
$i$th lattice point is $b(i)$ units above (resp. below) the terminal 
point of the $i$th step if the corresponding subscript is $N$ (resp. 
$Y$), 
as illustrated.
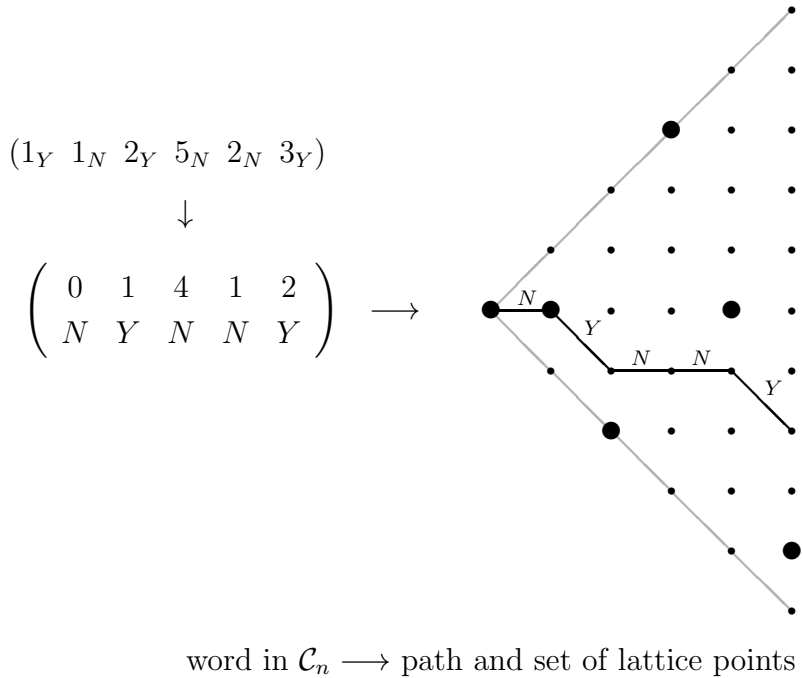
\Einheit=0.8cm
\[
\gray{ \Pfad(0,0),33333\endPfad
\Pfad(0,0),44444\endPfad }
\Pfad(0,0),14114\endPfad
\DuennPunkt(1,0)
\DuennPunkt(2,0)
\DuennPunkt(3,0)
\DuennPunkt(4,0)
\DuennPunkt(5,0)
\DuennPunkt(1,1)
\DuennPunkt(2,1)
\DuennPunkt(3,1)
\DuennPunkt(4,1)
\DuennPunkt(5,1)
\DuennPunkt(2,2)
\DuennPunkt(3,2)
\DuennPunkt(4,2)
\DuennPunkt(5,2)
\DuennPunkt(3,3)
\DuennPunkt(4,3)
\DuennPunkt(5,3)
\DuennPunkt(4,4)
\DuennPunkt(5,4)
\DuennPunkt(5,5)
\DuennPunkt(1,-1)
\DuennPunkt(2,-1)
\DuennPunkt(3,-1)
\DuennPunkt(4,-1)
\DuennPunkt(5,-1)
\DuennPunkt(2,-2)
\DuennPunkt(3,-2)
\DuennPunkt(4,-2)
\DuennPunkt(5,-2)
\DuennPunkt(3,-3)
\DuennPunkt(4,-3)
\DuennPunkt(5,-3)
\DuennPunkt(4,-4)
\DuennPunkt(5,-4)
\DuennPunkt(5,-5)
\DickPunkt(0,0)
\DickPunkt(1,0)
\DickPunkt(2,-2)
\DickPunkt(3,3)
\DickPunkt(4,0)
\DickPunkt(5,-4)
\Label\l{(1_{Y}\ \,1_{N}\ \,2_{Y}\ \,5_{N}\ \,2_{N}\ \,3_{Y})}(-5,2.6)
\Label\l{\downarrow}(-4.7,1.6)
\Label\l{\left(
\begin{array}{ccccc}
    0 & 1 & 4 & 1 & 2  \\
    N & Y & N & N & Y
\end{array}
\right) \ \ \longrightarrow
}(-4.1,0)
\Label\o{\textrm{{\scriptsize $N$}}}(0.6,-0.1)
\Label\o{\textrm{{\scriptsize $Y$}}}(1.7,-0.6)
\Label\o{\textrm{{\scriptsize $N$}}}(2.5,-1.1)
\Label\o{\textrm{{\scriptsize $N$}}}(3.5,-1.1)
\Label\o{\textrm{{\scriptsize $Y$}}}(4.7,-1.6)
\]
\centerline{word in $\C_{n} \longrightarrow $ path and set of lattice 
points}

The ordinates of the lattice points with a prepended 0 form a 
symmetric trapezoidal word (Section \ref{symTrap}). The path is 
redundant since it can be 
recovered from the lattice points and so we have an explicit bijection 
from $\C_{n}$ to symmetric trapezoidal words.

To code a height-labeled Dyck $n$-path by a word in $\C_{n}$, delete the last 
upstep $U$ with label 1 along with the last downstep and decrement labels 
after $U$ to get an $(n-1)$-path $P_{1}$ with a distinguished vertex $v$ 
(where the upstep 
was deleted). If $v$ is upstep-free (resp. upstep-incident) in 
$P_{1}$, the subscript on $w(n)$ is $Y$ (resp. $N$) and its value is 
the number of upstep-free (resp. upstep-incident) vertices weakly 
preceding $v$. Repeat on $P_{1}$ to get $w(n-1)$ and so on, ending 
with $w(1):=1_{Y}$. This gives a bijection that sends \#\,upstep-free 
vertices to \#\,subscripts $Y$.

Likewise, to code a Stirling $n$-permutation by a word in $\C_{n}$, delete 
the two $n$s (necessarily adjacent)
to get an $(n-1)$-permutation with a distinguished gap (possibly at 
either end). Recalling that the last gap is considered a descent, if 
the distinguished gap is a descent (resp. non-descent), the subscript on $w(n)$ is $Y$ (resp. $N$) 
and its value is 
the number of descent (resp. non-descent) gaps weakly 
preceding the distinguished gap. Repeat to get $w(n-1)$ and so on, 
again ending with $w(1):=1_{Y}$. This gives a bijection that sends 
\#\,descents to \#\,subscripts $Y$.

\noindent(5) This result is stated without proof in \cite{rior76}. 
A Stirling permutation can be built up in a unique way by starting 
with a plateau of two 1s, inserting a plateau 22 in one of the 3 gaps in 
--1--1--, then inserting a plateau 33 in one of 5 gaps and so on. 
Define a mapping from Stirling permutations of size $n$ to 
trapezoidal words $\mbf{w}=(w_{i})_{i=1}^{n}$ as follows. Set 
$w_{1}=1$. If 22 is placed between the two 1s, set $w_{2}=1$; if 22 is 
placed to the left of the ones, set $w_{2}=2$ (the smallest number 
not yet appearing in $\mbf{w}$), else $w_{2}=3$. In general, if $kk$ 
is placed inside a plateau, say the $i$th plateau (left to right), 
$w_{k}$ is the $i$th smallest number already appearing in $\mbf{w}$; 
otherwise $kk$ is placed in one of the remaining gaps, say in the 
$j$th of 
the remaining gaps, and $w_{k}$ is the $j$th smallest positive integer 
not yet appearing in $\mbf{w}$. For example, $5512234431 \rightarrow 
11442$ as follows.
\[
\begin{array}{ccccccccc}\\[-4mm]
    11 & \rightarrow & 1221 & \rightarrow & 122331 & \rightarrow & 12234431 & \rightarrow & 5512234431  \\
   w_{1}=1   &   & w_{2}=1 &  & w_{3}=4^{\,*} &  & w_{4}=4^{\,\dag} &  & w_{5}=2
\end{array}
\]
\hspace*{25mm} $^{*}$ 4 is the third smallest number not yet appearing in 
$\mbf{w}$ \\
\hspace*{25mm} $^{\dag}$ 4 is the second smallest number already appearing in 
$\mbf{w}$    

\vspace*{2mm}

\noindent It is easy to check that this algorithm defines a bijection 
from Stirling permutations of size $n$ to trapezoidal words of length 
$n$ and that it sends ``\#\:plateaus'' in the permutation to 
``\#\:distinct entries'' in the word.

\noindent(6) Translated to Stirling permutations using Janson's bijection 
described above, 
descents in a tree become strong descents in the permutation where a \emph{strong 
descent} is a descent $a>b$ involving the \emph{first} of the two occurrences 
of $b$ in the permutation. Note that the analogous notion of strong 
ascent is superfluous because all ascents in a Stirling permutation 
are strong. Consideration of the effect of inserting $nn$ into a 
Stirling permutation of size $n-1$ on the number of strong descents 
leads to the defining recurrence for the reversed second-order Eulerian triangle 
using the following fact, proved by induction: in a Stirling permutation 
of size $n$, \#\,strong descents + \#\,ascents $=n-1$.

\noindent(7) The number of upstep-free vertices in a Dyck $n$-path is related to the number of peaks: 
their sum is $n+1$. It follows that the number of height-labeled Dyck $n$-paths with $k$ peaks is 
$\SecondEulerianInline{n}{n+1-k}$, the second-order Eulerian triangle with 
reversed rows. 

\subsection{}\label{stirling}
\vspace*{-7mm}
\begin{equation}
\sum_{k=0}^{n} \StirlingCycle{n}{k}2^{n-k} = (2n-1)!!
    \label{firstascent1s}
\end{equation}
Here $\StirlingCycle{n}{k}$ is the Stirling cycle number (and 
$\StirlingCycle{0}{0}=1$).
This identity counts 

\noindent(1) Stirling permutations of size $n$ by number $k$ of 
left-to-right (LR) minima. (For 
example, $5\,5\,2\,3\,4\,4\,3\,2\,1\,1$ has 3 LR minima, namely 5,2,1), 

\noindent(2) increasing ordered trees of size $n$ by 
number $k$ of edges in the minimal path from the root as defined in Section 
\ref{stirliByFirst}, and

\noindent(3) height-labeled Dyck $n$-paths by number $k$ of 
upsteps in the first ascent with label 1.

\noindent \textbf{Proofs} 

\noindent(1) It suffices to exhibit a bijection from Stirling 
permutations $\si$
of size $n$ with $k$ left-to-right minima to pairs $(A,\tau)$ where 
$A$ is a subset of $[n-k]$ and $\tau$ is a permutation of $[n]$ with $k$ left-to-right minima
(recall the number of such $\tau$ is $\StirlingCycle{n}{k}$). Split 
$\si$ just before each LR minimum to write $\si$ as a concatenation 
$\si_{1},\si_{2},\ldots,\si_{k}$ of lists. Each $\si_{i}$ is itself a 
Stirling permutation and their support sets, taken in order, form a 
partition of $[n]$ into blocks whose smallest elements are decreasing left 
to right. In view of these observations, the general case follows from 
the special case $k=1$ by amalgamating subsets and concatenating 
permutations, both taken on the appropriate support sets. So suppose 
$k=1$, implying that $\si$ starts with a 1.
For 
$j\in[2,n]$, let $i_{j}$ denote the last entry preceding $j$ in $\si$ 
that is $<j$. Thus, for $\si=1\,2\,5\,5\,2\,1\,4\,4\,3\,3,\ 
i_{2},i_{3},i_{4}$ all $=1$ and $i_{5}=2$. Now take 
$A=\{j\in[2,n]\,:\,$both occurrences of $i_{j}$ in $\si$ precede the 
first occurrence of $j$\}, a subset of the $(n-1)$-element support set 
$[2,n]$. As for $\tau$, observe that $(i_{j})_{j=2}^{n} \in 
[1]\times[2]\times \cdots \times [n-1]$ and so corresponds to a permutation 
of the $(n-1)$-element set $[2,n]$. Prepend a 1 to this permutation to 
get the permutation $\tau$ of $[n]$ with $k=1$ LR minima. We leave the 
reader to verify that $\si$ can be recovered from the pair $(A,\tau)$.

\noindent(2) Adding a leaf $n$ to a tree of size $n-1$ preserves the 
length of the minimal path except when $n$ is added as a child of the 
leaf that terminates the minimal path. We thus have the recurrence
\[
u(n,k)=(2n-2)u(n-1,k)+u(n-1,k-1),
\]
satisfied by the summand because it reduces to the basic recurrence 
for the Stirling cycle numbers: 
$\StirlingCycle{n}{k}=(n-1)\StirlingCycle{n-1}{k} + 
\StirlingCycle{n-1}{k-1}$.

\noindent(3)  We 
give yet another bijection from  \HL Dyck paths to \IOTs. This one sends 
\#\,upsteps immediately preceding the $i$th downstep to the outdegree of 
vertex $i-1\ (1\le i \le n)$---in brief, full ascent sequence of 
path $\rightarrow$ fertility sequence of tree---and sends locations of 1s 
on first ascent to locations of LR minima among children of the root. 
In particular, \#\,1s = \#\,LR minima, and Janson's bijection identifies 
LR minima among children of the 
root with left-to-right minima in a Stirling permutation. So (3) will 
follow from (1). The bijection is an algorithm to generate the edge 
list. With $n=10$, and the path below as a working example,

\Einheit=0.6cm
\[
\Pfad(-10,0),33334434333444433444\endPfad
\SPfad(-10,0),11111111111111111111\endSPfad
\DuennPunkt(-10,0)
\DuennPunkt(-9,1)
\DuennPunkt(-8,2)
\DuennPunkt(-7,3)
\DuennPunkt(-6,4)
\DuennPunkt(-5,3)
\DuennPunkt(-4,2)
\DuennPunkt(-3,3)
\DuennPunkt(-2,2)
\DuennPunkt(-1,3)
\DuennPunkt(0,4)
\DuennPunkt(1,5)
\DuennPunkt(2,4)
\DuennPunkt(3,3)
\DuennPunkt(4,2)
\DuennPunkt(5,1)
\DuennPunkt(6,2)
\DuennPunkt(7,3)
\DuennPunkt(8,2)
\DuennPunkt(9,1)
\DuennPunkt(10,0)
\Label\o{\textrm{{\scriptsize 1}}}(-9.7,0.3)
\Label\o{\textrm{{\scriptsize 2}}}(-8.7,1.3)
\Label\o{\textrm{{\scriptsize 1}}}(-7.7,2.3)
\Label\o{\textrm{{\scriptsize 2}}}(-6.7,3.3)
\Label\o{\textrm{{\scriptsize 3}}}(-3.7,2.3)
\Label\o{\textrm{{\scriptsize 2}}}(-1.7,2.3)
\Label\o{\textrm{{\scriptsize 1}}}(-0.7,3.3)
\Label\o{\textrm{{\scriptsize 5}}}(0.3,4.3)
\Label\o{\textrm{{\scriptsize 1}}}(5.3,1.3)
\Label\o{\textrm{{\scriptsize 3}}}(6.3,2.3)
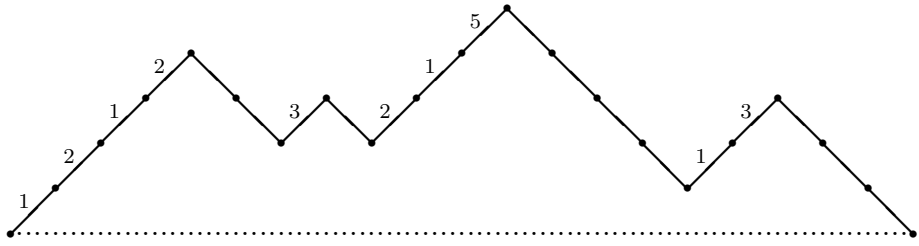
\Label\o{\textrm{{\small height-labeled Dyck 10-path}}}(0,-1.8)
\]

\vspace*{1mm}

\noindent list the total number of downsteps preceding the $ i$-th upstep 
for $ i=n,n-1,\ldots,1$ (column 2 in Fig. 3) and the label on the 
$ i$-th upstep (column 3). 

\vspace*{5mm}

$
\begin{array}{c|cccc}
   i & \#\,D\textrm{s} & \textrm{label on }
   & \textrm{candidate} & \textrm{selected}  \\[-2mm]
   &  <  i\textrm{-th } U  & i\textrm{-th } U  & \textrm{children}  & \textrm{edge}  \\
   \hline
    10 & 7 & 3 & 8\,9 \,10 & 7\,10  \\
    9 & 7 & 1 & 8\,9 & 7\,8  \\
    8 & 3 & 5 & 4\,5\,6\,7\,9 & 3\,9  \\
    7 & 3 & 1 & 4\,5\,6\,7 & 3\,4  \\
    6 & 3 & 2 & 5\,6\,7 & 3\,6  \\
    5 & 2 & 3 & 3\,5\,7 & 2\,7  \\
    4 & 0 & 2 & 1\,2\,3\,5 & 0\,2  \\
    3 & 0 & 1 & 1\,3\,5 & 0\,1  \\
    2 & 0 & 2 & 3\,5 & 0\,5  \\
    1 & 0 & 1 & 3 & 0\,3
\end{array}
$
\begin{pspicture}(-4,2)(4,4)
\psset{unit=.6cm}
\psline(-4,2)(-3,1)(0,0)(3,1)(3,2)(4,3)
\psline(-3,2)(-3,1)(-2,2)
\psline(-1,1)(0,0)(1,1)
\psline(2,3)(3,2)

\psdots(-4,2)(-3,1)(0,0)(3,1)(3,2)(4,3)(-3,2)(-2,2)(-1,1)(1,1)(2,3)

\rput(0,-0.4){\textrm{{\footnotesize $0$}}}
\rput(3.4,1){\textrm{{\footnotesize $2$}}}
\rput(-3.4,0.9){\textrm{{\footnotesize $3$}}}
\rput(-1,1.4){\textrm{{\footnotesize $5$}}}
\rput(1,1.4){\textrm{{\footnotesize $1$}}}
\rput(-4,2.4){\textrm{{\footnotesize $6$}}}
\rput(-3,2.4){\textrm{{\footnotesize $4$}}}
\rput(-2,2.4){\textrm{{\footnotesize $9$}}}
\rput(3.4,1.8){\textrm{{\footnotesize $7$}}}
\rput(2,3.4){\textrm{{\footnotesize $8$}}}
\rput(4,3.4){\textrm{{\footnotesize $10$}}}

\rput(-5.7,1.4){$\rightarrow$}

\end{pspicture}

\vspace*{4mm}

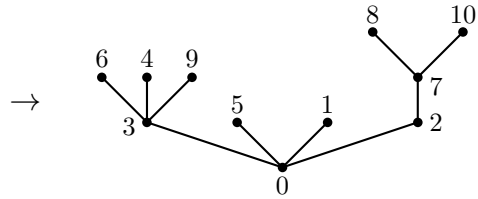
\centerline{Fig. 3}

\vspace*{4mm}

\noindent Each entry $ j $ in column 2 is the 
parent vertex of an edge. The candidates for its child are the entries 
of $ [j+1,n]$ not already having a parent. Initially $j=8$, no 
vertex has a parent, and 
the candidate children are 8,9,10. The corresponding label in column 3 
then determines the child, here label 3 selects the third candidate, 
namely 10, so 7\:10 becomes an edge; 10 now has a parent and is deleted 
from  the candidate set. The label 1 in the next row selects the 
first candidate, namely 8, and adds 7\:8 to the edge list. Proceed 
similarly to get all $n$ edges. Then the last column, 
read upwards, is the edge list in standard order, giving the tree 
shown.

The generating function 
$\sum_{k=1}^{n}\StirlingCycle{n}{k}2^{n-k}x^{n}$ for row $n$ is 
$\prod_{i=0}^{n-1}(x+2i)$, and the bivariate \gf $\sum_{n,k\ge 
0}\StirlingCycle{n}{k}2^{n-k}\frac{\textrm{{\,\normalsize 
$x$}}^{n}}{\textrm{{\small $n$}}!} y^{k}$ is 
$(1-2x)^{-y/2}$. See
\htmladdnormallink{A039683}{http://www.research.att.com:80/cgi-bin/access.cgi/as/njas/sequences/eisA.cgi?Anum=A102625}
 in OEIS.

\subsection{}
\vspace*{-5mm}
Consider the following modification of the minimal path from the root in 
\IOTs:
the right-then-minimal path is the path that goes from the 
root to its rightmost child, then follows minimal children to a leaf. 
Let $u(n,k)$ denote the number of \IOTs of size $n$ whose right-then-minimal 
path has length $k$, so that 
\[
\sum_{k=1}^{n}u(n,k)=(2n-1)!!.
\]
Then $u(n,k)$ satisfies the recurrence
\begin{eqnarray*}
    u(n,1) & = & n(2n-3)!!  \\
    u(n,k) & = & (2n-3)u(n-1,k)+u(n-1,k-1) \qquad 2\le k \le n,
\end{eqnarray*}
leading to a differential equation
for the \gf $F(x,y)=\sum_{n\ge k\ge 1}u(n,k)\frac{\textrm{{\,\normalsize $x$}}^{n}}{\textrm{{\small $n$}} 
\textrm{{\footnotesize !}} } y^{k}$:
\[
(1-2x)F_{x}(x,y)=\frac{y}{\sqrt{1-2x}} -(1-y)F(x,y),
\]
with solution
\[
F(x,y)=\frac{y\,\big(1-(1-2x)^{1-y/2}\big)}{(2-y)\sqrt{1-2x} }.
\]
The first few values of $u(n,k)$ are
\[
\begin{array}{c|ccccc}
	n^{\textstyle{\,\backslash \,k}} & 1 & 2 & 3 & 4 & 5 \\
\hline 
	1&    1 &   & & & \\
 	2&    2 & 1 & & &   \\
	3&    9 & 5  & 1 & &  \\
	4&    60 & 34  & 10 & 1 &   \\ 
	5&    525 & 298 & 104 & 17 & 1  \\

\end{array}
\]

\section{Refinements}\label{refine}

\vspace*{-5mm}
\subsection{}
\vspace*{-5mm}
The interpretations of Sections \ref{stirliByFirst} and \ref{stirling} have 
a common refinement
\[
\sum_{1\le j \le k \le 
n}2^{n-j}\StirlingCycle{k-1}{j-1}(n-1)^{\underline{n-k}}=(2n-1)!!.
\]
The summand is the number of \IOTs of size $n$ whose minimal path from 
the root ends at the leaf labeled $k$ and whose length is $j$: the 
number $u(n,k,j)$ of such paths satisfies the recurrence
\begin{eqnarray*}
    u(n,n,1) & = & (2n-2)!!,  \\
    u(n,k,1) & = & 0 \hspace*{44mm} k\ge 2, \\
    u(n,k,j) & = & (2n-2)u(n-1,k,j) \qquad 2\le j\le k \le n-1, 
\end{eqnarray*}
\begin{eqnarray*}
    u(n,n,j) & = & \sum_{i=j-1}^{n-1}u(n-1,i,j-1) \qquad 2 \le j \le n.
\end{eqnarray*}
The \gf $\sum_{1\le j \le k \le n}u(n,k,j)\frac{\textrm{{\normalsize $x$}}^{n-1}}{\textrm{{\small $(n-1)$}} 
\textrm{{\footnotesize !}} }y^{j-1}z^{k-1}$ is 
\[
F(x,y,z)=\frac{1}{(1-2x)(1-2xy)^{z/2}}.
\]

Likewise, the interpretations of Sections \ref{firstAscentLength} and \ref{stirling} have 
a common refinement:
\[
\sum_{1\le j \le k \le n} \StirlingCycle{k}{j} \binom{2n-k-1}{k-1}\,(2n-2k-1)!! = (2n-1)!!.
\]
The summand $v(n,k,j)$ is the number of height-labeled Dyck $n$-paths whose 
first ascent has length $k$ and contains $j$ 1s. The \gf $\sum_{0\le j 
\le k \le n} v(n,k,j) \frac{\textrm{{\normalsize $x$}}^{n}}{\textrm{{\small $n$}} 
\textrm{{\footnotesize !}} }y^{j}z^{k}$ is 
\[
\left(\frac{1-y-y\sqrt{1-2x}}{1-2y+2xy^{2}}\right)^{z}.
\]

\subsection{}
\vspace*{-5mm}
The result of Section \ref{firstAscentLength} can be refined to give 
the joint  distribution of first ascent 
length and first peak upstep label: let $u(n,j,k)$ denote the number 
of height-labeled Dyck $n$-paths 
whose first peak upstep has label $j$ and whose
initial ascent has length $k\ (n\ge k\ge j\ge 1)$.
Then, since $j$ is uniformly distributed over $[k]$,
$u(n,j,k)=(k-1)! \binom{2n-k-1}{k-1}\,(2n-2k-1)!!$ for $1\le j \le k$,
leading to the \gf $\sum_{n\ge 1,j\ge 1,k\ge 1}u(n,j,k)\frac{\textrm{{\normalsize $x$}}^{n-1}}{\textrm{{\small $(n-1)$}} 
\textrm{{\footnotesize !}} }y^{j-1}z^{k-1}=$
\[
\frac{1}{\sqrt{1-2x}\,\big(1-yz+yz\sqrt{1-2x}\,\big)\:\big(1-z+z\sqrt{1-2x}\,\big)}.
\]

\subsection{}
\vspace*{-5mm}
To refine the result of Section \ref{firstdescentsubsection} and deduce 
further identities, let 
$u(n,k,j)$ denote the number of height-labeled Dyck $n$-paths whose first ascent 
has length $=j$ and first descent has length $\ge k$. Thus 
$u(n,k)=\sum_{j=k}^{n}u(n,k,j)$. Since the first peak in such a path 
has $j$ possible labels, deleting this peak and its label shows that 
\begin{equation}
    u(n,k,j)=j u(n-1,k-1,j-1) \qquad \textrm{for $2\le k\le j \le n$}
    \label{firstascdes}
\end{equation}
In particular, $u(n,2,j)$ is $j$ times the number of size-$(n-1)$ height-labeled Dyck paths 
with first ascent of length $j-1$ and no restriction on the first descent. 
Hence, by the second interpretation of (\ref{firstascent}), 
$u(n,2,j)=j(j-1)(2n-2-j)^{\underline{j-2}}(2n-2j-1)!!$. This base 
case, together with (\ref{firstascdes}), yields
\[
u(n,k,j)=j^{\underline{k}}(2n-k-j)^{\underline{j-k}}(2n-2j-1)!!.
\]
Equating $\sum_{j=k}^{n}u(n,k,j)$ and $u(n,k)$ yields the identity
\[
    \sum_{j=k}^{n}j^{\underline{k}}(2n-k-j)^{\underline{j-k}}(2n-2j-1)!!=(k-1)!!(2n-k)!!.
\]
Two alternative forms of this identity, eliminating the 
double factorials, can be found by considering the cases where $k$ is even 
or odd separately:
\begin{equation}
\sum_{j=0}^{n}\binom{j+2m}{j}2^{j}\binom{2n-j}{n-j} = \binom{n+m}{n}4^{n}
    \label{ascdescentident1}
\end{equation}
with $k:=2m,\ n$ replaced by $n+2m$, and 
\begin{equation}
\sum_{j=0}^{n}\binom{j+2m+1}{j}2^{j+1}\binom{2n-j}{n-j} = 
\frac{m!(2n+2m+2)!}{n!(2m+1)!(n+m+1)!}
    \label{ascdescentident2}
\end{equation}
with $k:=2m+1,\ n$ replaced by $n+2m+1$.

This last identity is interesting because it provides another solution 
to Ira Gessel's 1987 Monthly Problem E3107:
\vspace*{-4mm}
\begin{quote}
    Show that $\frac{m!(2m+2n)!}{(2m)!n!(m+n)!}$ is an integer for 
    nonnegative integers $m,n$.
\end{quote}
Replace $m$ by $m-1$ in (\ref{ascdescentident2}), rearrange terms, and 
cancel a 2 to get
\[
\sum_{j=0}^{n}\binom{j+2m-1}{j}2^{j}\binom{2n-j}{n-j} = 
\frac{m!(2n+2m)!}{n!(2m)!(n+m)!},
\]
exhibiting Gessel's expression as a sum of integers, and this sum is 
different from that in A.\:A.\:Jagers' solution \cite{jagers}.
The case $m=0$ of (\ref{ascdescentident1}) has a simple combinatorial 
interpretation: it counts lattice paths of upsteps and downsteps of 
length $2n$ by number of ``returns to ground level''.


\begin{thebibliography}{99}

\bibitem{rior76} John Riordan, The blossoming of Schr\"{o}der's fourth problem, 
\emph{Acta Math.} \textbf{137} (1976), no. 1Ð2, 1-Ð16.


\bibitem{stirpoly} Ira Gessel and Richard Stanley, Stirling polynomials, \emph{J. Combinatorial Theory (A)} 
\textbf{24} (1978), 24--33.     
    
\bibitem{klazar1} Martin Klazar, 
Twelve countings with rooted plane trees, 
\emph{Eur. J. Combinatorics}, \textbf{18},
Issue 2 (1997), 195--210.Ê 

\bibitem{klazar2} Martin Klazar, 	
Addendum to ÒTwelve countings with rooted plane treesÓ,
\emph{Eur. J. Combinatorics} \textbf{18}, 
Issue 6 (1997), 739--740.Ê

\bibitem{janson08} Svante Janson, Plane recursive trees, Stirling permutations
and an urn model, \emph{Fifth Colloquium on Mathematics and Computer Science DMTCS Proceedings},
AI, (2008), 541-Ð548, \newline 
\htmladdnormallink{http://www.dmtcs.org/dmtcs-ojs/index.php/proceedings/article/view/dmAI0137}{http://www.dmtcs.org/dmtcs-ojs/index.php/proceedings/article/view/dmAI0137}.

\bibitem{ec2} Richard P.~Stanley, \emph{Enumerative Combinatorics} 
Vol.\,2, Cambridge University Press, 1999. Exercise 6.19 and related 
material on Catalan numbers are available online at
\htmladdnormallink{http://www-math.mit.edu/$\,\widetilde{\ }\,$rstan/ec/ }{http://www-math.mit.edu/~rstan/ec/}.


\bibitem{francon78} Jean Fran\c{c}on, 
Histoires de fichiers,\emph{ RAIRO Informat. ThŽor.} \textbf{12} (1978), no. 1, 49--62.

\bibitem{francon79} Jean Fran\c{c}on and G\'{e}rard Viennot, Permutations selon leurs pics, creux,
doubles mont\'{e}es et double descentes, nombres d'Euler et nombres de
Genocchi, \emph{Discrete Math.}  \textbf{28} (1979), no.  1, 21--35.

\bibitem{overhang}  Robert J. Marsh and Paul Martin, Tiling bijections 
between paths and Brauer diagrams, preprint, 4 Jun 2009,
\htmladdnormallink{http://arxiv.org/abs/0906.0912v1}{http://arxiv.org/abs/0906.0912v1}.

\bibitem{oeis}
The \htmladdnormallink{On-Line Encyclopedia of Integer 
Sequences}{http://www.research.att.com:80/~njas/sequences/},
founded and maintained by Neil J.~Sloane.

\bibitem{jagers} A.\,A.\,Jagers, solution for E3107, \emph{Amer. Math. Monthly},\textbf{ 95}, 
No. 1 (Jan., 1988), 53--54.


\end{thebibliography}
\end{document}